\documentclass{article}%
\usepackage{amsfonts}
\usepackage{graphicx}
\usepackage{amsmath}
\usepackage{amssymb}%
\providecommand{\U}[1]{\protect\rule{.1in}{.1in}}
\newtheorem {theorem}{Theorem}[section]
\newtheorem {proposition}{Proposition}[section]
\newtheorem {corollary}{Corollary}[section]

\newtheorem{example}{Example}[section]
\newtheorem{lemma}{Lemma}[section]

\newtheorem{remark}{Remark}[section]

\begin{document}

\begin{center}
{\LARGE Exact Moderate and Large Deviations for Linear Processes}

\bigskip

\centerline{\today}

\bigskip Magda Peligrad$^{a}$, Hailin Sang$^{b}$, Yunda Zhong$^{c}$ and Wei
Biao Wu$^{c}$

\bigskip$^{a}$ Department of Mathematical Sciences, University of Cincinnati,
PO Box 210025, Cincinnati, OH 45221-0025, USA. E-mail address: peligrm@ucmail.uc.edu

\bigskip$^{b}$ Department of Mathematics, University of Mississippi,
University, MS 38677-7071, USA. E-mail address: sang@olemiss.edu

\bigskip$^{c}$ Department of Statistics, University of Chicago, Chicago, IL
60637, USA. E-mail address: ydzhong@galton.uchicago.edu,
wbwu@galton.uchicago.edu \bigskip
\end{center}

\bigskip\textbf{Abbreviated Title: }{\Large Exact Deviations for Linear
Processes}

\begin{center}
\bigskip

\textbf{Abstract}
\end{center}

Large and moderate deviation probabilities play an important role in many
applied areas, such as insurance and risk analysis. This paper studies the
exact moderate and large deviation asymptotics in non-logarithmic form for
linear processes with independent innovations. The linear processes we analyze
are general and therefore they include the long memory case. We give an
asymptotic representation for probability of the tail of the normalized sums
and specify the zones in which it can be approximated either by a standard
normal distribution or by the marginal distribution of the innovation process.
The results are then applied to regression estimates, moving averages,
fractionally integrated processes, linear processes with regularly varying
exponents and functions of linear processes. We also consider the computation
of value at risk and expected shortfall, fundamental quantities in risk theory
and finance.

\bigskip


\footnotetext[3]{This work was supported in part by the Taft Research Center
at the University in Cincinnati. In addition, Magda Peligrad was supported in
part by the NSA grant H98230-11-1-0135 and NSF DMS-1208237. Wei Biao Wu by NSF
grants DMS-0906073 and DMS-1106790.}

\section{Introduction and notations}

\bigskip

Let $(\xi_{i})_{i\in\mathbb{Z}}$ be a sequence of independent and identically
distributed centered random variables with finite second moment and $c_{n i}$
a sequence of constants. This paper focuses on the moderate and large
deviations in non-logarithmic form for the linear process of the form
\begin{equation}
S_{n}=\sum_{i=1}^{k_{n}}c_{ni}\xi_{i}. \label{gln}%
\end{equation}
This class of linear processes is versatile enough to help analyzing
regression estimates, moving averages that include long memory processes,
linear processes with regularly varying coefficients and fractionally
integrated processes.

Our goal is to find an asymptotic representation for the tail probabilities of
the normalized sums defined by (\ref{gln}). Estimations of deviation
probabilities occur in a natural way in many applied areas including insurance
and risk analysis.

Specifically, we aim to find a function $N_{n}(x)$ such that, as
$n\rightarrow\infty$,
\begin{equation}
{\frac{{\mathbb{P}(S_{n}\geq x\sigma_{n})}}{{{N_{n}(x)}}}}%
=1+o(1),\mbox{ where }\sigma_{n}^{2}=\Vert S_{n}\Vert_{2}^{2}=\mathbb{E}%
\xi_{1}^{2}\sum_{i=1}^{k_{n}}c_{ni}^{2}. \label{eq:nonlgd}%
\end{equation}
If $x\geq0$ is fixed, then (\ref{eq:nonlgd}) becomes the well-known central
limit theorem by letting $N_{n}(x)=1-\Phi(x)$, where $\Phi(x)$ is the standard
normal distribution function. In this paper we call $\mathbb{P}(S_{n}
/\sigma_{n}\geq x)$ the \textit{moderate} or \textit{large deviation}
probabilities depending on the speed of convergence $x=x_{n}\rightarrow\infty
$. These tail probabilities of rare events can be very small. Here we call
(\ref{eq:nonlgd}) the \textit{exact approximation}, which is more accurate
than the logarithmic version
\begin{equation}
{\frac{{\log\mathbb{P}(S_{n}/\sigma_{n}\geq x)}}{{{\log N_{n}(x)}}}}=1+o(1),
\label{eq:lgd}%
\end{equation}
which is often used in the literature in the context of large or moderate
deviation. For example, suppose $\mathbb{P}(S_{n} / \sigma_{n}\geq x)
=10^{-4}$ and $N_{n}(x)=10^{-5}$; then their logarithmic ratio is $0.8$, which
does not appear to be very different from $1$, while the ratio for the exact
version (\ref{eq:nonlgd}) is as big as $10$. A multiplicative factor of this
order can cause substantially different industrial standards in designing
projects that can survive natural disasters. The logarithmic version
(\ref{eq:lgd}) is incapable of effectively characterizing the differences
between the tail probabilities.

As early as 1929, Khinchin considered the problem of moderate and large
deviation probabilities in non-logarithmic form for independent Bernoulli
random variables. The first large deviation probability result appeared in S.
Nagaev (1965). A. Nagaev (1969) studied large deviation probabilities of
i.i.d. random variables with regularly varying tails. Mikosch and A. Nagaev
(1998) applied the large deviation probabilities for heavy-tailed random
variables to insurance mathematics. The review work on this topic can be found
in S. Nagaev (1979) and Rozovski (1993). Rubin and Sethuraman (1965),
Slastnikov (1978) and Frolov (2005) considered the moderate or large
deviations for arrays of independent random variables. S. Nagaev (1979)
presented the following very useful result: in (\ref{gln}) assume $k_{n}=n$,
$c_{ni}\equiv1$, and that $\xi_{i}$ has a regularly varying right tail. i.e.
\begin{equation}
\mathbb{P}(\xi_{0}\geq x)=\frac{h(x)}{x^{t}}\text{ as }x\rightarrow
\infty\mbox{ for some }t>2, \label{tail1}%
\end{equation}
where $h(x)$ is a slowly varying function (Bingham, Goldie and Teugels, 1987).
Namely, $\lim_{x\rightarrow\infty}h(\lambda x)/h(x)=1$ for all $\lambda>0$. If
in addition, for some $p>2$, $\xi_{0}$ has absolute moment of order $p$, then
\begin{equation}
\mathbb{P}(\sum_{i=1}^{n}\xi_{i}\geq x\sigma_{n})=(1-\Phi
(x))(1+o(1))+n\mathbb{P}(\xi_{0}\geq x\sigma_{n})(1+o(1)) \label{eq:NagaevJ16}%
\end{equation}
for $n\rightarrow\infty$ and $x\geq1$. Note that (\ref{eq:NagaevJ16}) implies
(\ref{eq:nonlgd}) with
\begin{equation}
N_{n}(x)=(1-\Phi(x))+n\mathbb{P}(\xi_{0}\geq x\sigma_{n}).
\label{eq:NagaevJ21}%
\end{equation}
Hence if $1-\Phi(x)=o[n\mathbb{P}(\xi_{0}\geq x\sigma_{n})]$ (resp.
$n\mathbb{P}(\xi_{0}\geq x\sigma_{n})=o(1-\Phi(x))$), then in (\ref{eq:nonlgd}%
) we can also choose $N_{n}(x)=1-\Phi(x)$ (resp. $N_{n}(x)=n\mathbb{P}(\xi
_{0}\geq x\sigma_{n})$).

The study of moderate and large deviation probabilities in non-logarithmic
form for dependent random variables is still in its initial stage. Ghosh
(1974) considered moderate deviations for $m$-dependent random variables. Chen
(2001) obtained a moderate deviation result for Markov processes. Grama (1997)
and Grama and Haeusler (2006) investigated the martingale case. Wu and Zhao
(2008) studied moderate deviations for stationary processes which applies to
many time series models. However the result in the latter paper can only be
applied to linear processes with short memory and their transformations.

For analyzing linear processes with long memory and for obtaining other
interesting applications, we study processes of type (\ref{gln}). Under mild
conditions on the coefficients, we shall point out the zones in which the
deviation probabilities can be approximated either by a standard normal
distribution or by using the distribution of $\xi_{0}$. Our main result is
that (\ref{eq:NagaevJ16}) holds in our case with
\[
N_{n}(x)=(1-\Phi(x))+\sum_{i=1}^{k_{n}}\mathbb{P}(c_{ni}\xi_{0}\geq
x\sigma_{n}).
\]
The paper has the following structure. Section \ref{sec:mainrlt} presents a
general moderate and large deviation result and various applications. Section
\ref{sec:num} illustrates the results of a numerical study. In Section
\ref{sec:proof} we prove the results. In the Appendix we give some auxiliary
results and we also mention some known facts needed for the proofs.

Before stating our results we introduce the notations that will be used
throughout this paper: $a_{n}\sim b_{n}$ means that $\lim_{n\rightarrow\infty
}a_{n}/b_{n}=1$, $a_{n}=O(b_{n})$ and also $a_{n}\ll b_{n}$ means
$\limsup_{n\rightarrow\infty} a_{n}/b_{n}<\infty$; $a_{n}=o(b_{n})$ if
$\lim_{n \rightarrow\infty} a_{n}/b_{n}=0$. By $\Vert X\Vert_{p}$ we denote
$(\mathbb{E}|X|^{p})^{1/p}.$ The notation $l(\cdot)$, $h(\cdot)$ and
$\ell(\cdot)$ denote slowly varying functions. By convention $0/0$ is
interpreted as $0$.

\section{Main Results}

\label{sec:mainrlt}

\bigskip

Throughout the paper, we assume that:

\textbf{Condition A.} $(\xi_{i})_{i\in\mathbb{Z}}$, are i.i.d. centered random
variables with finite second moment, $\sigma^{2}=\mathbb{E}\xi_{0}^{2}$.

\bigskip

\subsection{General linear processes}

\bigskip

Our first results apply to general linear processes of type (\ref{gln}) with
i.i.d. innovations. For $c_{ni}>0$ and $t>0$, we define
\begin{equation}
B_{nt}= \sum_{i=1}^{k_{n}}c_{ni}^{t}, \label{defAnt}%
\end{equation}%
\begin{equation}
\sigma_{n}^{2}=var(S_{n})=B_{n2}\mathbb{E}\xi_{0}^{2}, \label{def-sigma}%
\end{equation}
 and
\begin{equation}
D_{nt}=B_{n2}^{-t/2}B_{nt}. \label{def-Dn}%
\end{equation}
The basic assumption in all our results is the uniform asymptotic
negligibility of the variance of individual summands, namely
\begin{equation}
\max_{1\leq i\leq k_{n}}c_{ni}^{2}/\sigma_{n}^{2}\rightarrow0. \label{cni}%
\end{equation}
Our first theorem extends Nagaev's result in (\ref{eq:NagaevJ16}) to general
linear processes.

\begin{theorem}
\label{mix} Assume that $(\xi_{i})_{i\in\mathbb{Z}}$ satisfies Condition A,
and for a certain $t>2$ it satisfies the right tail condition (\ref{tail1}).
Moreover, for a certain $p>2$, $\Vert\xi_{0}\Vert_{p}<\infty.$ Assume also
that $c_{ni}>0$ and (\ref{cni}) is satisfied. Let $(x_{n})_{n\ge1}$ be any
sequence such that for some $c>0$ we have $x_{n}\geq c$ for all $n$. Then, as
$n\rightarrow\infty$,
\begin{equation}
\mathbb{P}\left(  S_{n}\geq x_{n}\sigma_{n}\right)  =(1+o(1))\sum_{i=1}%
^{k_{n}}\mathbb{P}(c_{ni}\xi_{0}\geq x_{n}\sigma_{n})+(1-\Phi(x_{n}))(1+o(1)).
\label{MD+LD}%
\end{equation}

\end{theorem}

\begin{remark}
To be precise, in relation (\ref{MD+LD}) as well as in (\ref{LD}) and
(\ref{MD}) below, by $o(1)$ we understand a function which depends on $x_{n}$
and on the underlying distribution, with the property that
its limit as $n\rightarrow\infty$ is zero. Each $o(1)$ may represent a
different function. The sequence $(x_{n})_{n\geq1}$ may be bounded or may
converge to infinity.

\end{remark}

\begin{corollary}
\label{Remark 1} Under the conditions of Theorem \ref{mix} for $x_{n}\geq
a(\ln D_{nt}^{-1})^{1/2}$ with $a>2^{1/2}$ we have%
\begin{equation}
\mathbb{P}(S_{n}\geq x_{n}\sigma_{n})=(1+o(1))\sum_{i=1}^{k_{n}}%
\mathbb{P}(c_{ni}\xi_{0}\geq x_{n}\sigma_{n})\text{ as }n\rightarrow\infty.
\label{LD}%
\end{equation}
On the other hand, if $0<x_{n}\leq b(\ln D_{nt}^{-1})^{1/2}$ with $b<2^{1/2},$
we have
\begin{equation}
\mathbb{P}\left(  S_{n}\geq x_{n}\sigma_{n}\right)  =(1-\Phi(x_{n}%
))(1+o(1))\text{ as }n\rightarrow\infty. \label{MD}%
\end{equation}

\end{corollary}

\begin{remark}
Notice that (\ref{LD}) and (\ref{MD}) assert different approximations for the
tail probability $\mathbb{P}(S_{n}\geq x\sigma_{n})$: moderate behavior for
$x=x_{n}$ smaller than a threshold, when we can approximate this probability
by using a normal distribution. On the other hand we have a large deviation
type of behavior for $x$ larger than another threshold. The behavior at the
boundary $\sqrt2 (\ln D_{nt}^{-1})^{1/2}$ is more subtle and it depends on the
slowly varying function $h(\cdot)$. For the special case in which $\lim_{x
\to\infty} h(x) \to h_{0}>0$, we have
\begin{equation}
{\frac{{\mathbb{P}(S_{n}\geq x\sigma_{n})} }{{N_{n}(x)}}} = 1+o(1),
\mbox{ where } N_{n}(x) = (1-\Phi(x)) + {\frac{{h_{0}}}{{(\sigma x)^{t}}}}D_{nt}.
\label{eq:NM1911}%
\end{equation}
If $x \geq a(\ln D_{nt}^{-1})^{1/2}$ with $a>2^{1/2},$ then $N_{n}(x) \sim
h_{0} D_{nt} /(\sigma  x)^{t}$.
\end{remark}

The proofs of these results are based on a separate study of the behaviors of
type (\ref{LD}) or (\ref{MD}), which is of independent interest. As a matter
of fact, we shall see in the next two theorems that a result similar to
(\ref{LD}) holds without the assumption of the finite moment of order $p>2$
while the moderate deviation (\ref{MD}) does not require a regularly varying right tail.

\begin{theorem}
\label{LargeD1} Assume that $(\xi_{i})_{i\in\mathbb{Z}}$ satisfies Condition
A, and for a certain $t>2$ it satisfies (\ref{tail1}). Let $c_{ni}>0$ be a
sequence of constants satisfying (\ref{cni}). Then, for any sequence
$x_{n}\geq C_{t}(\ln D_{nt}^{-1})^{1/2}$ with $C_{t}>e^{t/2}(t+2)/\sqrt{2}$
the large deviation result (\ref{LD}) holds.
\end{theorem}

As a counterpart to this result we shall formulate now the moderate deviation
bound. 
\begin{theorem}
\label{ModerateD} Assume that $(\xi_{i})_{i\in\mathbb{Z}}$ satisfies Condition
A and for a certain $p>2$, $\Vert\xi_{0}\Vert_{p}<\infty$. Assume that
(\ref{cni}) is satisfied. If $x_{n}^{2}\leq2\ln(D_{np}^{-1})$ then the
moderate deviation result (\ref{MD}) holds.
\end{theorem}

\subsection{Applications to linear regression estimates}

Many statistical procedures, such as estimation of regression coefficients,
produce linear statistics of type (\ref{gln}). See for instance Chapter 9 in
Beran (1994), for the case of parametric regression, or the paper by Robinson
(1997), where kernel estimators are used for nonparametric regression. Here we
consider the simple parametric regression model $Y_{i}=\beta\alpha_{i}+\xi
_{i}$, where $\xi_{i}$ are i.i.d. centered errors with $\mathbb{E}\xi_{1}%
^{2}=\sigma^{2}$, $(\alpha_{i})$ is
a sequence of positive real numbers and $\beta$ is the parameter of interest.
The least squares estimator $\hat{\beta}_{n}$ of $\beta$, based on a sample of
size $n,$ satisfies%

\begin{equation}
S_{n}:=\hat{\beta}_{n}-\beta=\frac{1}{\sum_{i=1}^{n}\alpha_{i}^{2}}\sum
_{i=1}^{n}\alpha_{i}\xi_{i}, \label{eq:J22235}%
\end{equation}
so, the representation of type (\ref{gln}) holds with $c_{ni}=\alpha_{i}%
/(\sum_{i=1}^{n}\alpha_{i}^{2})$. Denote $A_{nt}=\sum_{i=1}^{n}\alpha_{i}^{t}%
$. Notice that $var(S_{n})=\sigma^{2}/A_{n2}.$ Assume
\begin{equation}
\lim_{n\rightarrow\infty}A_{n2}^{-1}\max_{1\leq i\leq n}\alpha_{i}^{2}=0.
\label{as1}%
\end{equation}
As an immediate consequence of Theorem \ref{mix}, we obtain: 
\begin{corollary}
\label{cor:mdldlf} (i) Assume that $(\xi_{i})_{i\in\mathbb{Z}}$ and $x=x_{n}$
satisfies the conditions in Theorem \ref{mix}. Under assumption (\ref{as1}),
we have
\begin{align*}
\mathbb{P}(\hat{\beta}_{n}-\beta &  \geq x\sigma/A_{n2}^{1/2})=\\
(1+o(1))\sum_{i=1}^{n}\mathbb{P}(\xi_{i}  &  \geq x\sigma A_{n2}^{1/2}%
/\alpha_{i})+(1+o(1))(1-\Phi(x)).
\end{align*}
(ii) If $x>0$ and $x^{2}\leq2\ln(A_{n2}^{t/2}/A_{nt})$, under the conditions
in Theorem \ref{mix}, we have
\[
\mathbb{P}(\hat{\beta}_{n}-\beta\geq x\sigma/A_{n2}^{1/2})=(1+o(1))(1-\Phi
(x)).
\]
(iii) If $x>0$ and $x^{2}\geq C_{t}^{2}\ln(A_{n2}^{t/2}/A_{nt})$ with
$C_{t}^{2}>2$, under the conditions in Theorem \ref{mix}, then
\[
\mathbb{P}(\hat{\beta}_{n}-\beta\geq x\sigma/A_{n2}^{1/2})=(1+o(1))\sum
_{i=1}^{n}\mathbb{P}(\xi_{i}\geq x\sigma A_{n2}^{1/2}/\alpha_{i}).
\]

\end{corollary}

Similar results as in Theorems \ref{LargeD1} and \ref{ModerateD} can also be
easily formulated.

Theorems \ref{mix}, \ref{LargeD1} and \ref{ModerateD} are also applicable to
the nonlinear regression model $y_{i}=g(x_{i})+\xi_{i}$, $1\leq i\leq n$,
where $g(x)$ is an unknown function and $\xi_{i}$ is the noise. Let $x_{i}$ be
the deterministic design points. Then the Nadaraya-Watson estimate $\hat
{g}_{n}$ satisfies
\[
\hat{g}_{n}(x)-\mathbb{E}\hat{g}_{n}(x)=\sum_{i=1}^{n}c_{ni}(x)\xi_{i}%
\]
where, letting $K$ be a kernel function and $h_{n}$ be bandwidths, the
weights
\[
c_{ni}(x)=K\left(  \frac{x_{i}-x}{h_{n}}\right)  /\sum_{i=1}^{n}K\left(
\frac{x_{i}-x}{h_{n}}\right)  .
\]
Therefore it is of the type (\ref{gln}).

\bigskip

\subsection{Application to moving averages}

\label{sec:alpJ16} \label{sec:alpJ16 copy(1)}

\bigskip

We now consider the sum $S_{n}=\sum_{k=1}^{n}X_{k}$, where
\begin{equation}
X_{k}=\sum_{j=-\infty}^{\infty}a_{k-j}\xi_{j}. \label{ln}%
\end{equation}
We assume that $\sum_{i\in{\mathbb{Z}}}a_{i}^{2}<\infty,$ which is the
necessary and sufficient condition for the existence of $X_{1}$. Observe that
$S_{n}=\sum_{i=-\infty}^{\infty}b_{ni}\xi_{i}$ is of form (\ref{gln}) with
\begin{equation}
b_{ni}=a_{1-i}+\cdots+a_{n-i} \label{def-cni}%
\end{equation}
and $k_{n}=\infty$. Assume $b_{ni}>0$ for all $i$ and let 
\begin{equation}
U_{nt}=(\sum_{i}b_{ni}^{2})^{-t/2}\sum_{i}b_{ni}^{t}. \label{def U}%
\end{equation}
In the corollary below, this quantity will replace $D_{nt}$ from definition
(\ref{def-Dn}) and $b_{ni}$ will replace the $c_{ni}$ in Subsection 2.1. Define $\sigma_n^2=\mathbb{E}\xi_0^2 \sum_{i}b_{ni}^{2}$. We
know from Peligrad and Utev (1997) that under the assumption $\sigma_{n}%
^{2}\rightarrow\infty$ we have
\begin{equation}
\sigma_{n}^{-2}\sup_{i}b_{ni}^{2}\rightarrow0\text{ as }n\rightarrow\infty.
\label{bni}%
\end{equation}
Therefore condition (\ref{cni}) is automatically satisfied. As a corollary of
Theorems \ref{mix}, \ref{LargeD1} and \ref{ModerateD} we obtain:

\begin{corollary}
\label{LinearLDMD} Assume that $(X_{n})_{n\geq1}$ is defined by (\ref{ln}) and
$\sigma_{n}^{2}\rightarrow\infty$. 
\newline(i) Assume that $(\xi_{i}%
)_{i\in\mathbb{Z}}$ and $x_n$ satisfy the conditions of Theorem \ref{mix} and
$b_{ni}>0$. Then (\ref{MD+LD}) holds. Corollary \ref{Remark 1} is also valid
for the partial sum of (\ref{ln}). 
\newline(ii) Let $(\xi_{i})_{i\in
\mathbb{Z}}$ be as in Theorem \ref{LargeD1}. Assume $b_{ni}>0$. Then the large
deviation result (\ref{LD}) holds for the sequence $x_{n}\geq C_{t}(\ln
U_{nt}^{-1})^{1/2}$ with $C_{t}>e^{t/2}(t+2)/\sqrt{2}$. 
\newline(iii) Assume
$(\xi_{i})_{i\in\mathbb{Z}}$ is as in Theorem \ref{ModerateD}.
Then the moderate deviation result (\ref{MD}) holds for $x_{n}^{2}\leq
2\ln(U_{np}^{-1})$.
\end{corollary}

Note that this corollary applies to general linear processes including the
long memory processes with $\sum_{i}|a_{i}|=\infty.$ Asymptotic properties for
long memory processes can be quite different from those of processes with
short memory, partially because the variance of the partial sum goes to
infinity at an order different than $n$; see for example, Ho and Hsing (1997),
Robinson (2003), Doukhan, Oppenheim and Taqqu (2003) among others. Hall (1992)
gave a Berry-Esseen bound for the convergence rate in the central limit theorem.

We shall apply now this corollary to the important particular case of causal
long-memory processes with
\begin{equation}
a_{i}=l(i+1)(1+i)^{-r},\text{ }i\geq0,\text{ with }1/2<r<1,\text{ and }%
a_{i}=0\text{ in rest}. \label{eq:rvaiJ22252}%
\end{equation}
Here $l(\cdot)$ is a slowly varying function where the results can be given in
a more precise form. Notice that in this particular case
\begin{equation}
X_{k}=\sum_{j=-\infty}^{k}a_{k-j}\xi_{j}. \label{defx}%
\end{equation}
Let $a_{0}=1$. This case of long memory linear processes covers the well-known
fractional ARIMA processes (cf. Granger and Joyeux; 1980, Hosking, 1981),
which play an important role in financial time series modeling and
application. As a special case, let $0<d<1/2$ and $B$ be the backward shift
operator with $B\varepsilon_{k}=\varepsilon_{k-1}$ and consider
\[
X_{k}=(1-B)^{-d}\xi_{k}=\sum_{i\geq0}a_{i}\xi_{k-i},\text{ where }a_{i}%
=\frac{\Gamma(i+d)}{\Gamma(d)\Gamma(i+1)}.
\]
For this example we have $\lim_{n\rightarrow\infty}a_{n}/n^{d-1}=1/\Gamma(d)$.
Note that these processes have long memory because $\sum_{j\geq0}%
|a_{j}|=\infty.$

\begin{corollary}
\label{RegularLDMD} Assume (\ref{eq:rvaiJ22252}). If $(\xi_{i})_{i\in
\mathbb{Z}}$ satisfies the conditions of \ Theorem \ref{mix} then
(\ref{MD+LD}) holds. In particular (\ref{LD}) holds for $x_n\geq c_{1}(\ln
n)^{1/2}$ with $c_{1}>(t-2)^{1/2}$ while (\ref{MD}) holds, provided $0<x_n\leq
c_{2}(\ln n)^{1/2}$ with $c_{2}<(t-2)^{1/2}$.
\end{corollary}

For this case Theorems \ref{LargeD1} and \ref{ModerateD} give:

\begin{corollary}
\label{Separate}(i) Let $(\xi_{i})_{i\in\mathbb{Z}}$ be as in Theorem
\ref{LargeD1}. Then (\ref{LD}) holds for $x_n>c_{1}(\ln n)^{1/2}$ with
$c_{1}>(t-2)^{1/2}e^{t/2}(t+2)/2$.\newline(ii) Let $(\xi_{i})_{i\in\mathbb{Z}%
}$ be as in Theorem \ref{ModerateD}. Then (\ref{MD}) holds if $x_n^{2}%
\leq(p-2)\ln n$.
\end{corollary}

\subsection{Application to risk measures}

\bigskip

In risk theory and finance, value at risk (VaR) and expected shortfall (ES)
play a fundamental role; see Jorion (2006), Holton (2003), McNeil et al
(2005), Acerbi and Tasche (2002) among others. Mathematically, they are
equivalent to quantiles and tail conditional expectations. In practice one is
most interested in their extremal behavior which corresponds to tail
quantiles. Despite their importance, however, their computation can be quite
difficult and the related asymptotic justification is far from being trivial.

Here we shall apply Theorem \ref{mix} and provide approximate formulae for
extremal quantiles and tail conditional expectations for $S_{n}$ defined by
(\ref{gln}). Under the assumption $\lim_{x\rightarrow\infty}h(x)\rightarrow
h_{0}>0$, by (\ref{eq:NM1911}) and Theorem \ref{mix},
\[
\mathbb{P}\left(  S_{n}\geq x\sigma_{n}\right)  =(1+o(1)){\frac{{h_{0}}%
}{(\sigma{x)^{t}}}}D_{nt}+(1-\Phi(x))(1+o(1)).
\]
Given the tail probability $\alpha\in(0,1)$, let $q_{\alpha,n}$ be the upper
$\alpha$-th quantile of $S_{n}$. Namely $\mathbb{P}(S_{n}\geq q_{\alpha
,n})=\alpha$. Elementary calculations show that $q_{\alpha,n}$ can be
approximated by $x_{\alpha}\sigma_{n}$ in the sense that $\lim_{n\rightarrow
\infty}x_{\alpha}\sigma_{n}/q_{\alpha,n}=1$, where $x=x_{\alpha}$ is the
solution to the equation
\[
{\frac{{h_{0}}}{(\sigma{x)^{t}}}}D_{nt}+(1-\Phi(x))=\alpha.
\]
In particular,  if $\alpha\leq h_{0}D_{nt}((a\sigma)^{2}\ln D_{nt}^{-1})^{-t/2}$ with
$a>2^{1/2}$, then, by Corollary \ref{Remark 1}, we can approximate
$q_{\alpha,n}$ by $\sigma^{-1}(h_{0}D_{nt}/\alpha)^{1/t}\sigma_{n}=\sigma
^{-1}{(B_{nt}h_{0}/\alpha)^{1/t}.}$ The approximation is understood in the
sense that $\sigma^{-1}{(B_{nt}h_{0}/\alpha)^{1/t}}/q_{\alpha,n}\rightarrow1$
as $n\rightarrow\infty$, and the tail conditional expectation or expected
shortfall is computed as
\begin{align*}
\mathbb{E}(S_{n}|S_{n} &  \geq q_{\alpha,n})=\frac{q_{\alpha,n}\mathbb{P}%
(S_{n}\geq q_{\alpha,n})+{\int_{q_{\alpha,n}}^{\infty}\mathbb{P}(S_{n}\geq
w)dw}}{{\mathbb{P}(S_{n}\geq q_{\alpha,n})}}\\
&  \sim q_{\alpha,n}+\frac{{q_{\alpha,n}}}{{t-1}}={\frac{t{q_{\alpha,n}}%
}{{t-1}}}\sim\sigma^{-1}B_{nt}^{1/t}{\frac{t{(h_{0}/\alpha)^{1/t}}}{{t-1}}}.
\end{align*}
We emphasize that, without the exact moderate deviation principle in  Corollary
 \ref{Remark 1}, the validity of the above equivalence cannot be guaranteed. To the
best of our knowledge, our example might be the only case that one can obtain
explicit asymptotic expressions for VaR and ES for sums of dependent random variables.

\subsection{Functionals of linear processes}

\bigskip

In this subsection we shall use the result from (ii) of Corollary
\ref{Separate} to study the moderate deviation for nonlinear transformations
of linear processes. Let $K$ be a transformation which is measurable and
$\mathbb{E}K(X_{0})=0$. Let
\[
H_{n}=\sum_{i=1}^{n}K(X_{i})\text{ where }X_{i}\text{ is defined by
(\ref{defx})}.
\]
For example, if $K(X_{0})=I(X_{0}\leq\tau)-\mathbb{P}(X_{0}\leq\tau)$, then
$H_{n}/n$ becomes the empirical process. If $X_{i}$ is short memory, namely
$a_{i}$ are absolutely summable, then we can apply the moderate deviation
principle in Wu and Zhao (2008). However, the result in the latter paper is
not applicable for long-range dependent processes. Despite its importance in
risk analysis, the problem of moderate deviation under strong dependence has
been rarely studied in the literature.

Here we shall establish such a principle in the context of nonlinear
transforms of linear processes. First, we introduce some necessary notation
for this section. Let $\mathcal{F}_{n} = (\cdots, \xi_{n-1},\xi_{n})$ be the
shift process and define the projection operator $\mathcal{P}_{i} \cdot=
\mathbb{E}(\cdot|\mathcal{F}_{i})-\mathbb{E}(\cdot|\mathcal{F}_{i-1})$. Denote
the truncated processes $X_{n,k} = \mathbb{E}(X_{n}|\mathcal{F}_{k})$. Now
define the functions $K_{n}(w)=\mathbb{E}[K(w+X_{n}-X_{n,0})]$ and $K_{\infty
}(w)=\mathbb{E}[K(w+X_{n})]$. We consider transformations $K$ with
$\kappa:=K_{\infty}^{\prime}(0)\neq0$. Define
\[
S_{n,1}=\sum_{i=1}^{n}[K(X_{i})-\kappa X_{i}]=H_{n}-\kappa S_{n}%
,\mbox{ where }S_{n}=\sum_{i=1}^{n}X_{i}.
\]
Then $H_{n}=\kappa S_{n}+S_{n,1}$. For a function $g$, let $g(w;\lambda
)=\sup_{|y|\leq\lambda}|g(w+y)|$ be the local maximal function. Denote the
collection of functions with second order partial derivatives by
$\mathbb{C}^{2}(\mathbb{R})$. We need the following regularity condition.

\bigskip

\textbf{Condition B.} Let $2\leq q<p\leq2q$ and assume $\Vert\xi_{0}\Vert
_{p}<\infty$. Assume $K_{n}\in C^{2}(\mathbb{R})$ for all large $n$ and that
for some $\lambda>0$,
\[
\sum_{i=0}^{2}\Vert K_{n-1}^{(i)}(X_{n,0};\lambda)\Vert_{q}+\Vert|\xi
_{1}|^{p/q}K_{n-1}(X_{n,1})\Vert_{q}+\Vert\xi_{1}K_{n-1}^{\prime}%
(X_{n,1})\Vert_{q}=O(1).
\]
A version of Condition B with $q=2$ is used in Wu (2006). We shall establish
the following moderate deviation result. For $1/2<r<1$ and $1/2\leq v<1$
define%
\begin{gather*}
\chi(v,r)=v\max(r-r/v,\,1/2-r,\,r-1),\\
\omega(r)=\mathrm{argmin}_{1/2\leq v<1}\chi(v,r)\text{ and }\rho
(r)=-\chi(\omega(r),r).
\end{gather*}

\begin{theorem}
\label{thmfunc} Assume that Condition B holds with $q=p\omega(r)$ and the
conditions of Corollary \ref{Separate} (ii) are satisfied. Let $c$ be such
that $0<c\leq p-2$ and $c<2p\rho(r).$ Then if $x\leq c\ln n$, we have
\begin{equation}
\mathbb{P}(H_{n}\geq|\kappa|\sigma_{n}x)=(1-\Phi(x))(1+o(1))\text{ as
}n\rightarrow\infty. \label{MDH}%
\end{equation}

\end{theorem}

\begin{remark}
\label{var} As mentioned in the proof of Theorem \ref{thmfunc} in Section
\ref{sec:prf24}, (\ref{MDH}) is still valid if the normalizing constant
$|\kappa|\sigma_{n}$ therein is replaced by $\sqrt{var(H_{n})}.$
\end{remark}

\begin{remark}
Theorem \ref{thmfunc} only asserts a moderate deviation with the Gaussian
range. It is unclear whether the approximation of type (\ref{LD}) holds. We
pose it as an open problem.
\end{remark}

\begin{remark}
An explicit form for $\omega(r)$ can be obtained. If $r\geq3/4$, then
$\omega(r)=r$. If $r<3/4$, then $\omega(r)=r/(2r-1/2)$. If $2p\rho(r)\geq
p-2$, then the moderate deviation in (\ref{MDH}) has the same range as for
$S_{n}$. The latter happens, for example, if $r=3/4$ and $2<p<16/5,$ since in
this case $2p\rho(3/4)\geq p-2$.
\end{remark}

\begin{example}
As an application to empirical processes, let $K(X)=I(X\leq\tau)-\mathbb{P}%
(X\leq\tau)$, where $\tau\in\mathbb{R}$ is fixed. Let $X_{n}=\xi_{n}%
+\sum_{i=1}^{\infty}a_{i}\xi_{n-i}=:\xi_{n}+Y_{n-1}$, where $\Vert\xi_{0}%
\Vert_{p}<\infty$, $p>2$, and its density function $f_{\xi}$ satisfies
\begin{equation}
\sup_{u}[f_{\xi}(u)+|f_{\xi}^{\prime}(u)|]<\infty. \label{eq:f1f2}%
\end{equation}
Then $K_{1}(w)=F_{\xi}(\tau-w)-F_{X}(\tau)$, where $F_{\xi}$ is the
distribution function of $\xi_{i}$. Under (\ref{eq:f1f2}), we clearly have
$\sup_{w}[|K_{1}^{\prime}(w)|+|K_{1}^{\prime\prime}(w)|]<\infty$. Observe that
we have the identity: for $n\geq1$,
\[
K_{n}(w)=\mathbb{E}K_{1}(w+a_{1}\xi_{n-1}+a_{2}\xi_{n-2}+\ldots+a_{n-1}\xi
_{1}).
\]
Hence $\sup_{n}\sup_{w}[|K_{n}^{\prime}(w)|+|K_{n}^{\prime\prime}(w)|]<\infty
$. So Condition B holds for any $\lambda$ since $\xi_{n}\in L^{p}$, $p>2$.
\end{example}

\bigskip

\section{A Numerical Study}

\label{sec:num}

\bigskip

In this section we shall design a numerical study of the accuracy of the large
deviation (\ref{LD}), normal approximation (\ref{MD}) and also the estimate
(\ref{MD+LD}). In particular, we shall study the accuracy of the
approximations in Corollary \ref{RegularLDMD}. In general it is very
time-consuming to calculate tail probabilities by Monte-Carlo simulation,
especially if they are small. One may need to carry out astronomically large
amount of computations to obtain reasonably well approximations.

Here we shall approach the problem from a different angle. We let $X_{j}%
=\sum_{i=1}^{\infty}a_{i}\xi_{j-i}$, where $\xi_{i}$, $i\in\mathbb{Z}$, have
Student's t-distribution with degree of freedom $\nu=3$, and $a_{i}=i^{-0.9}$.
Let $S_{n}=\sum_{i=1}^{n}X_{i}$ with $n=300$. Note that the characteristic
function of $\xi_{i}$ is
\begin{equation}
\varphi(t)={\frac{{(\sqrt{\nu}|t|)^{\nu/2}K_{\nu/2}(\sqrt{\nu}|t|)}}%
{{\Gamma(\nu/2)2^{\nu/2-1}}}}, \label{inversion}%
\end{equation}
where $K_{\nu/2}$ is the Bessel function (see Hurst (1995)). Then the
characteristic function of $S_{n}$ is
\[
\varphi_{S_{n}}(t)=\prod_{j\in\mathbb{Z}}\varphi(b_{nj}t)
\]
and by the inversion formula,
\[
\mathbb{P}(S_{n}\leq x)-\mathbb{P}(S_{n}\leq x^{\prime})={\frac{1}{{2\pi}}%
}\int_{-\infty}^{\infty}{\frac{{e^{\sqrt{-1}yx}-e^{\sqrt{-1}yx^{\prime}}}%
}{{\sqrt{-1}y}}}\varphi_{S_{n}}(y)dy.
\]
In the above equation let $x^{\prime}=0$. Since $\xi_{j}$ is symmetric,
$\mathbb{P}(S_{n}\leq0)=1/2$. In our numerical study we shall use
(\ref{inversion}) to compute the probability $\mathbb{P}(S_{n}>x)$.

\begin{figure}[h]
\label{fig1}
\par
\begin{center}
\includegraphics[width=8cm]{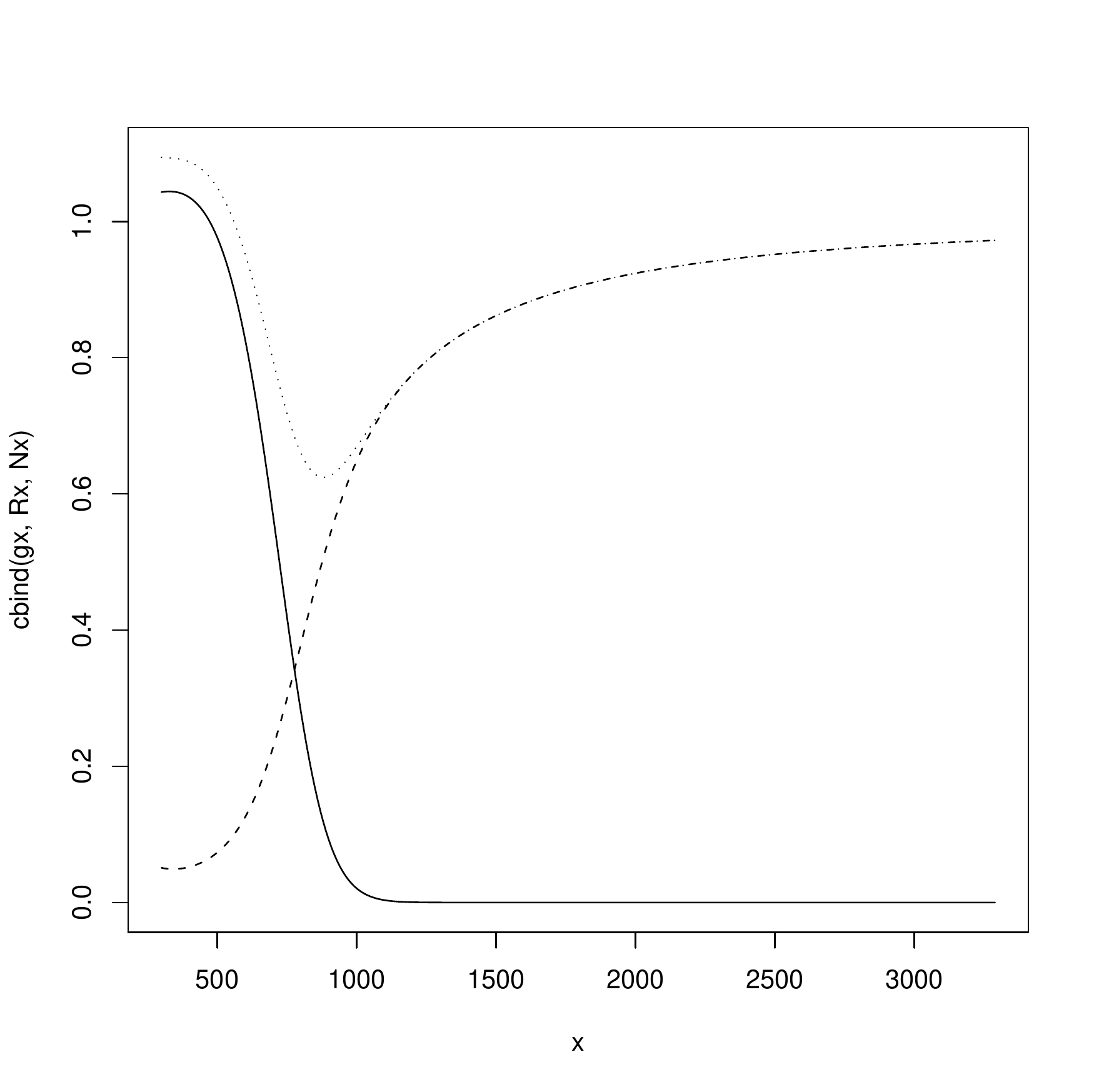}
\end{center}
\par
\textbf{Fig. 1.} Tail approximation $R(x)$ (dashed curve), Gaussian
approximation $g(x)$ (solid curve) and their sum (dotted curve) for
long-memory processes with Student $t(3)$ innovations.\end{figure}In Figure
$1$ we report the ratios $R(x):=\sum_{i}\mathbb{P}(b_{ni}\xi_{0}\geq x) /
\mathbb{P}(S_{n}>x)$ and $g(x):=(1-\Phi(x/\sigma_{n})) / \mathbb{P}(S_{n}>x)$;
see (\ref{LD}) with $c_{ni}=b_{ni}$. We can interpret $R(x)$ (resp. $g(x)$) as
tail (resp. Gaussian) approximation. As expected from Corollary
\ref{RegularLDMD}, the Gaussian approximation is better if $x$ is small, while
the tail probability $R(x)$ approximation is better when $x$ is big. In the
intermediate region we approximate by their sum.

\bigskip

\section{Proofs}

\label{sec:proof}

\subsection{Preliminary approximations}

Let $(X_{i})_{1\leq i\leq n}$ be independent random variables. We shall
approximate the tail distribution of partial sums by the tail of the sums of
truncated random variables and a term involving the tail probabilities of
individual summands. We use the following notations:
\[
S_{n}=\sum_{i=1}^{n}X_{i}\text{, }S(j)=\sum_{i\neq j}^{n}X_{i}%
\]
and for $x>0$ and $\varepsilon>0$ we set
\begin{equation}
X_{i}^{(\varepsilon x)}=X_{i}I(X_{i}<\varepsilon x)\text{, }S_{n}%
^{^{(\varepsilon x)}}=\sum_{i=1}^{n}X_{i}^{(\varepsilon x)}\text{ and }%
S_{n}^{^{(\varepsilon x)}}(j)=\sum_{i\neq j}^{n}X_{i}^{(\varepsilon x)}.
\label{def trunc}%
\end{equation}
We shall prove the following key lemma that will be further exploited to
approximate the tail distribution of $\mathbb{P}(S_{n}\geq x)$ in terms of the
sum of the truncated random variables and the tail distributions of the
individual summands.

\begin{lemma}
\label{approx} For any $0<\eta<1,$ and $\varepsilon>0$ such that
$1-\eta>\varepsilon$ we have%
\begin{gather*}
|\mathbb{P}(S_{n}\geq x)-\mathbb{P}(S_{n}^{(\varepsilon x)}\geq x)-\sum
_{j=1}^{n}\mathbb{P}(X_{j}\geq(1-\eta)x)|\leq\\
4{\Large (}\sum_{j=1}^{n}\mathbb{P}(X_{j}\geq\varepsilon x){\Large )}%
^{2}+3\sum_{j=1}^{n}\mathbb{P}(X_{j}\geq\varepsilon x)(\mathbb{P}%
(|S_{n}(j)|>\eta x)\\
+\sum_{j=1}^{n}\mathbb{P}((1-\eta)x\leq X_{j}<(1+\eta)x).
\end{gather*}

\end{lemma}

\textbf{Proof.} We decompose the event $\{S_{n}\geq x\}$ according to
$\max_{i\neq j}X_{i}<\varepsilon x$ or $\max_{i\neq j}X_{i} \geq\varepsilon
x,$ and the last one can happen if exactly one of the variables is larger than
$\varepsilon x$ or at least two variables exceed $\varepsilon x.$ Formally,%

\begin{gather*}
\mathbb{P}(S_{n}\geq x)=\sum_{j=1}^{n}\mathbb{P}(S_{n}\geq x,\text{ }X_{j}%
\geq\varepsilon x,\text{ }\max_{i\neq j}X_{i}<\varepsilon x)\\
+\mathbb{P}(%
{\displaystyle\bigcup\limits_{1\leq i\leq n-1}}\;\;
{\displaystyle\bigcup\limits_{i+1\leq j\leq n}}
\{S_{n}\geq x,\text{ }X_{j}\geq\varepsilon x,\text{ }X_{i}\geq\varepsilon
x\})\\
+\mathbb{P}(S_{n}\geq x,\text{ }\max_{1\leq i\leq n}X_{i}<\varepsilon
x)=A+B+C=\sum_{j=1}^{n}A_{j}+B+C.
\end{gather*}
The term $B$ can be easily majorated by
\[
B\leq\sum_{i=1}^{n-1}\sum_{j=i+1}^{n}\mathbb{P}(X_{j}\geq\varepsilon
x)\mathbb{P(}X_{i}\geq\varepsilon x)\leq{\Large (}\sum_{j=1}^{n}%
\mathbb{P}(X_{j}\geq\varepsilon x){\Large )}^{2}.
\]
We analyze now the first term. We introduce a new parameter $\eta>0$. Since
for any two events $A$ and $B$ we have $|P(A)-P(B)|\leq P(AB^{\prime
})+P(A^{\prime}B),$ (here the prime stays for the complement), for each $j$ we
have
\begin{align*}
|A_{j}-\mathbb{P}(X_{j}  &  \geq(1-\eta)x)|\leq\mathbb{P}(S_{n}\geq x,\text{
}X_{j}\geq\varepsilon x,\text{ }X_{j}<(1-\eta)x)\\
+\mathbb{P}(X_{j}  &  \geq(1-\eta)x,\text{ }S_{n}<x)+\mathbb{P}(X_{j}%
\geq(1-\eta)x,\text{ }X_{j}<\varepsilon x)\\
+\mathbb{P}(X_{j}  &  \geq(1-\eta)x,\text{ }\max_{i\neq j}X_{i}\geq\varepsilon
x)=I+II+III+IV.
\end{align*}
We treat each term separately. By independence and since $S_{n}\geq x$ and
$X_{j}<(1-\eta)x$ imply $S_{n}(j)\geq\eta x,$ we derive
\[
I\leq\mathbb{P}(X_{j}\geq\varepsilon x)\mathbb{P}(S_{n}(j)\geq\eta x).\text{ }%
\]
The second term is treated in the following way:%
\begin{align*}
II  &  \leq\mathbb{P}((1-\eta)x\leq X_{j}<(1+\eta)x)+\mathbb{P}(X_{j}%
\geq(1+\eta)x,\text{ }S_{n}<x)\\
&  \leq\mathbb{P}((1-\eta)x\leq X_{j}<(1+\eta)x)+\mathbb{P}(X_{j}\geq
(1+\eta)x)\mathbb{P}(-S_{n}(j)\geq\eta x).
\end{align*}
Since $1-\eta>\varepsilon$ the third term is: $III=0$. By independence, the
forth term is
\[
IV=\mathbb{P}(X_{j}\geq(1-\eta)x)\mathbb{P}(\max_{i\neq j}X_{i}\geq\varepsilon
x).
\]
Overall, by the previous estimates and because $1-\eta>\varepsilon$, we obtain%
\begin{gather*}
|A-\sum_{j=1}^{n}\mathbb{P}(X_{j}\geq(1-\eta)x)|\leq2\sum_{j=1}^{n}%
\mathbb{P}(X_{j}\geq\varepsilon x)(\mathbb{P}(|S_{n}(j)|>\eta x)\\
+{\Large (}\sum_{j=1}^{n}\mathbb{P}(X_{j}\geq\varepsilon x){\Large )}^{2}%
+\sum_{j=1}^{n}\mathbb{P}((1-\eta)x\leq X_{j}<(1+\eta)x).
\end{gather*}
It remains to analyze the last term, $C$. Notice that
\begin{align*}
|C-\mathbb{P}(S_{n}^{(\varepsilon x)}  &  \geq x)|=\mathbb{P}(S_{n}%
^{(\varepsilon x)}\geq x)-\mathbb{P}(S_{n}^{(\varepsilon x)}\geq x,\text{
}\max_{1\leq i\leq n}X_{i}<\varepsilon x)\\
&  =\mathbb{P}(S_{n}^{(\varepsilon x)}\geq x,\max_{1\leq i\leq n}X_{i}%
\geq\varepsilon x).
\end{align*}
Now we treat this term by the same arguments we have already used, by dividing
the maximum in two parts:
\begin{gather*}
\mathbb{P}(S_{n}^{(\varepsilon x)}\geq x,\max_{1\leq i\leq n}X_{i}%
\geq\varepsilon x)=\sum_{j=1}^{n}\mathbb{P}(S_{n}^{(\varepsilon x)}\geq
x,\text{ }X_{j}\geq\varepsilon x,\text{ }\max_{i\neq j}X_{i}<\varepsilon x)\\
+\mathbb{P}(%
{\displaystyle\bigcup\limits_{1\leq i\leq n-1}}\;\;
{\displaystyle\bigcup\limits_{i+1\leq j\leq n}}
\{S_{n}^{(\varepsilon x)}\geq x,\text{ }X_{j}\geq\varepsilon x,\text{ }%
X_{ni}\geq\varepsilon x\})=\sum_{j=1}^{n}F_{j}+G.
\end{gather*}
The last term, $G$ is majorated exactly as $B$. As for the first term, we
notice that because $X_{j}\geq\varepsilon x$ the term $X_{j}^{(\varepsilon
x)}$ does not appear in the sum, and by independence we obtain
\begin{align*}
F_{j}  &  =\mathbb{P}(S_{n}^{(\varepsilon x)}(j)\geq x,\text{ }X_{j}%
\geq\varepsilon x,\text{ }\max_{i\neq j}X_{i}<\varepsilon x)\\
&  \leq\mathbb{P}(S_{n}^{(\varepsilon x)}(j)\geq x)\mathbb{P(}X_{j}%
\geq\varepsilon x).
\end{align*}
Now, clearly we have%
\begin{align*}
\mathbb{P}(S_{n}^{(\varepsilon x)}(j)  &  \geq x)\leq\mathbb{P}(\max_{i}%
X_{i}\geq\varepsilon x)+\mathbb{P}(S_{n}^{(\varepsilon x)}(j)\geq x,\text{
}\max_{i}X_{i}<\varepsilon x)\\
&  =\mathbb{P}(\max_{i}X_{i}\geq\varepsilon x)+\mathbb{P}(S_{n}(j)\geq
x,\text{ }\max_{i}X_{i}<\varepsilon x),
\end{align*}
implying that
\[
\sum_{j=1}^{n}F_{j}\leq\sum_{j=1}^{n}\mathbb{P(}X_{nj}\geq\varepsilon
x)(\mathbb{P}(\max_{i}X_{i}\geq\varepsilon x)+\mathbb{P}(S_{n}(j)\geq x)).
\]
Overall,%
\[
|C-\mathbb{P}(S_{n}^{(\varepsilon x)}\geq x)|\leq2{\Large (}\sum_{j=1}%
^{n}\mathbb{P}(X_{j}\geq\varepsilon x){\Large )}^{2}+\sum_{j=1}^{n}%
\mathbb{P(}X_{j}\geq\varepsilon x)\mathbb{P}(S_{n}(j)\geq x).
\]
By gathering all the information above and taking into account that
\begin{gather*}
|\mathbb{P}(S_{n}\geq x)-\mathbb{P}(S_{n}^{(\varepsilon x)}\geq x)-\sum
_{j=1}^{n}\mathbb{P}(X_{j}\geq(1-\eta)x)|\leq\\
|A-\sum_{j=1}^{n}\mathbb{P}(X_{j}\geq(1-\eta)x)|+|C-\mathbb{P}(S_{n}%
^{(\varepsilon x)}\geq x)|+|B|,
\end{gather*}
the lemma is established. $\diamondsuit${}\textbf{ }

\bigskip The following similar lemma is for the sum of infinite many terms.

\begin{lemma}
\label{approxinfty} Let $1-\eta>\varepsilon>0$ and $x>0$; let $X_{1}%
,X_{2},\cdots,$ be independent random variables. Assume that the sum
$S=\sum_{i=1}^{\infty}X_{i}$ exists almost surely. Let $S_{(j)}=S-X_{j},\;X_{i}%
^{(\varepsilon x)}=X_{i}I(X_{i}<\varepsilon x).$ Then$\;S^{(\varepsilon
x)}=\sum_{i=1}^{\infty}X_{i}^{(\varepsilon x)}$ exists almost surely and
\begin{gather*}
|\mathbb{P}(S\geq x)-\mathbb{P}(S^{(\varepsilon x)}\geq x)-\sum_{j=1}^{\infty
}\mathbb{P}(X_{j}\geq(1-\eta)x)|\leq\\
4{\Large (}\sum_{j=1}^{\infty}\mathbb{P}(X_{j}\geq\varepsilon x){\Large )}%
^{2}+3\sum_{j=1}^{\infty}\mathbb{P}(X_{j}\geq\varepsilon x)(\mathbb{P}%
(|S(j)|>\eta x)\\
+\sum_{j=1}^{\infty}\mathbb{P}((1-\eta)x\leq X_{j}<(1+\eta)x).
\end{gather*}

\end{lemma}

\textbf{Proof. }By Kolmogorov's three-series theorem, $S^{(\varepsilon
x)}=\sum_{i=1}^{\infty}X_{i}^{(\varepsilon x)}$ converges almost surely. Let
$\Omega_{0}\in\Omega$ with $\mathbb{P}(\Omega_{0})=1$ be the set that both
$\sum_{i=1}^{\infty}X_{i}$ and $\sum_{i=1}^{\infty}X_{i}^{(\varepsilon x)}$
converge. Hence on $\Omega_{0}$, we understand $S(\omega)$ as just the sum
$\sum_{i=1}^{\infty}X_{i}(\omega)$. Then following the proof of Lemma
\ref{approx}, we have Lemma \ref{approxinfty}. $\diamondsuit${}\textbf{ }
\bigskip

If $S_{n}$ is stochastically bounded, i.e., $\lim_{K\rightarrow\infty}\sup
_{n}\mathbb{P}(|S_{n}|>K)=0,$ the approximation in Lemma \ref{approx} has a
simple asymptotic form.

\begin{proposition}
\label{stockB}Assume that $S_{n}$ is stochastically bounded, the variables are
centered and $x_{n}\rightarrow\infty$. Then for any $0<\eta<1,$ and
$\varepsilon>0$ such that $1-\eta>\varepsilon,$ we have%
\begin{gather}
|\mathbb{P}(S_{n}\geq x_{n})-\mathbb{P}(S_{n}^{(\varepsilon x_{n})}\geq
x_{n})-\sum_{j=1}^{n}\mathbb{P}(X_{j}\geq(1-\eta)x_{n})|\leq\label{ineqprop}\\
o(1)\sum_{j=1}^{n}\mathbb{P}(X_{j}\geq\varepsilon x_{n})+\sum_{j=1}%
^{n}\mathbb{P}((1-\eta)x_{n}\leq X_{j}<(1+\eta)x_{n}),\nonumber
\end{gather}
where $o(1)$ depends on the sequence $x_{n},$ $\eta$ and $\varepsilon$ and
converges to $0$ as $n\rightarrow\infty.$
\end{proposition}

\textbf{Proof. } We just notice that for independent centered random
variables, if $S_{n}$ is stochastically bounded, by L\'{e}vy inequality
(Inequality 1.1.3 in de la Pe\~{n}a and Gin\'{e} 1999), we have $\max_{1\leq
i\leq n}|X_{i}|$ is stochastically bounded too. By taking into account that
$|S_{n}(j)|\leq|S_{n}|+\max_{1\leq i\leq n}|X_{i}|,$ and using the fact that
$x_{n}\rightarrow\infty$ as $n\rightarrow\infty$ we obtain
\begin{gather*}
\sum_{j=1}^{n}\mathbb{P}(X_{j}\geq\varepsilon x_{n})\mathbb{P}(|S_{n}%
(j)|\geq\eta x_{n})\leq\max_{1\leq j\leq n}\mathbb{P}(|S_{n}(j)|\geq\eta
x_{n})\sum_{j=1}^{n}\mathbb{P}(X_{j}\geq\varepsilon x_{n})\\
\leq\left(  \mathbb{P}(|S_{n}|\geq\eta x_{n}/2)+\mathbb{P(}\max_{1\leq i\leq
n}|X_{i}|\geq\eta x_{n}/2)\right)  \sum_{j=1}^{n}\mathbb{P}(X_{j}%
\geq\varepsilon x_{n})\\
=o(1)\sum_{j=1}^{n}\mathbb{P}(X_{j}\geq\varepsilon x_{n})\text{ as
}n\rightarrow\infty.
\end{gather*}
Then, by independence
\begin{align*}
\mathbb{P}(\max_{1\leq j\leq n}|X_{j}|  &  \geq\varepsilon x_{n})\\
&  =\mathbb{P}(|X_{1}|\geq\varepsilon x_{n})+\sum_{k=2}^{n}\mathbb{P}%
(\max_{1\leq j\leq k-1}|X_{j}|<\varepsilon x_{n})\mathbb{P(}|X_{k}%
|\geq\varepsilon x_{n})\\
&  \geq\mathbb{P}(\max_{1\leq j\leq n}|X_{j}|<\varepsilon x_{n})\sum_{k=1}%
^{n}\mathbb{P}(|X_{j}|\geq\varepsilon x_{n}),
\end{align*}
which gives
\begin{align*}
{\Large (}\sum_{j=1}^{n}\mathbb{P}(|X_{j}|  &  \geq\varepsilon x_{n}%
){\Large )}^{2}\leq\frac{\mathbb{P}(\max_{1\leq j\leq n}|X_{j}|\geq\varepsilon
x_{n})}{\mathbb{P}(\max_{1\leq j\leq n}|X_{j}|<\varepsilon x_{n})}\sum
_{j=1}^{n}\mathbb{P}(|X_{j}|\geq\varepsilon x_{n})\\
&  =o(1)\sum_{j=1}^{n}\mathbb{P}(|X_{j}|\geq\varepsilon x_{n})\text{ as
}n\rightarrow\infty,
\end{align*}
since $x_{n}\rightarrow\infty$ as $n\rightarrow\infty$ and $\max_{1\leq j\leq
n}|X_{j}|$ is stochastically bounded. $\diamondsuit${}

\begin{remark}
\label{stochBremark} Based on Lemma \ref{approxinfty}, it is easy to verify
that Proposition \ref{stockB} is still valid if we extend the sums up to infinity.
\end{remark}

\subsection{Proof of Theorem \ref{LargeD1}}

It is convenient to normalize by the variance of partial sum and we shall
consider without restricting the generality that%

\begin{equation}
\mathbb{E}\xi_{0}^{2}=1,\;\;\sum_{i=1}^{k_{n}}c_{ni}^{2}=1\text{ and }%
\max_{1\leq i\leq k_{n}}c_{ni}^{2}\rightarrow0. \label{sum1}%
\end{equation}
Then we have $\sum_{i=1}^{k_{n}}c_{ni}^{t}\leq\max_{1\leq i\leq k_{n}}%
c_{ni}^{t-2}\rightarrow0$ implying that $D_{nt}^{-1}\rightarrow\infty.$
Moreover, the sequence $\sum_{i=1}^{k_{n}}c_{ni}\xi_{i}$ is stochastically
bounded and we analyze the two terms of the right side and the last term of
the left side in Proposition \ref{stockB}. Let $x_{n}\rightarrow\infty$ as
$n\rightarrow\infty.$ In order to ease the notation we shall denote $x=x_{n},$
but we keep in mind that $x$ depends on $n$ and tends to infinite with $n$. By
taking into account that $x/c_{ni}\geq x\rightarrow\infty$ and $h$ is a slowly
varying function we notice first that for any $a>0$
\[
\lim_{x\rightarrow\infty}\max_{1\leq i\leq k_{n}}|\frac{h(ax/c_{ni}%
)}{h(x/c_{ni})}-1|=0.
\]
We derive for any $|\gamma|<1$ fixed
\begin{gather*}
|\sum_{i=1}^{k_{n}}c_{ni}^{t}(h(\frac{x}{c_{ni}})-h((1+\gamma)\frac{x}{c_{ni}%
}))|\leq\\
\sum_{i=1}^{k_{n}}c_{ni}^{t}h(\frac{x}{c_{ni}})|1-\frac{h((1+\gamma)x/c_{ni}%
)}{h(x/c_{ni})}|=o(1)\sum_{i=1}^{k_{n}}c_{ni}^{t}h(\frac{x}{c_{ni}}),\text{ as
}n\rightarrow\infty,
\end{gather*}
implying that
\begin{align*}
\frac{\sum_{i=1}^{k_{n}}\mathbb{P}(c_{ni}\xi_{i}\geq(1\pm\eta)x)}{\sum
_{i=1}^{k_{n}}\mathbb{P}(c_{ni}\xi_{i}\geq x)}  &  =\frac{\sum_{i=1}^{k_{n}%
}c_{ni}^{t}h((1\pm\eta)x/c_{ni})}{(1\pm\eta)^{t}\sum_{i=1}^{k_{n}}c_{ni}%
^{t}h(x/c_{ni})}\rightarrow1\text{ }\\
\text{when }n  &  \rightarrow\infty\text{ followed by }\eta\rightarrow0.
\end{align*}
Then, we also have
\[
\frac{\sum_{i=1}^{k_{n}}\mathbb{P}((1-\eta)x\leq c_{ni}\xi_{i}<(1+\eta
)x)}{\sum_{i=1}^{k_{n}}\mathbb{P}(c_{ni}\xi_{i}\geq x)}\rightarrow0\text{ as
}n\rightarrow\infty\text{ and }\eta\rightarrow0.
\]
Similarly, for every $\varepsilon>0$ fixed we have that
\[
\frac{\sum_{i=1}^{k_{n}}\mathbb{P}(c_{ni}\xi_{i}\geq\varepsilon x)}{\sum
_{i=1}^{k_{n}}\mathbb{P}(c_{ni}\xi_{i}\geq x)}=\frac{\sum_{i=1}^{k_{n}}%
c_{ni}^{t}h(\varepsilon x/c_{ni})}{\varepsilon^{t}\sum_{i=1}^{k_{n}}c_{ni}%
^{t}h(x/c_{ni})}\rightarrow\frac{1}{\varepsilon^{t}}\text{ as }n\rightarrow
\infty,
\]
 and then,%
\[
\sum_{i=1}^{k_{n}}\mathbb{P}(c_{ni}\xi_{i}\geq\varepsilon x)\ll\sum
_{i=1}^{k_{n}}\mathbb{P}(c_{ni}\xi_{i}\geq x)\text{ as }n\rightarrow\infty.
\]
So far, for any $\varepsilon>0$ fixed, by letting $n\rightarrow\infty$ first
and after that, passing with $\eta$ to $0$, we deduce by the above
consideration combined with Proposition \ref{stockB}~that
\begin{equation}
\mathbb{P}(S_{n}\geq x)=\sum_{i=1}^{k_{n}}\mathbb{P}(c_{ni}\xi_{i}\geq
x)(1+o(1))+\mathbb{P}(S_{n}^{(\varepsilon x)}\geq x)\text{ }\ \text{as
}n\rightarrow\infty. \label{ineqLMD}%
\end{equation}
It remains to study the term $\mathbb{P}(S_{n}^{(\varepsilon x)}\geq x).$ We
shall base this part of the proof on Corollary 1.7 in S. Nagaev (1979), given
in the Appendix, which we apply with $m>t,$ that will be selected later.
Because we assume $\mathbb{E}(\xi_{0}^{2})=1$ and $\sum_{i=1}^{k_{n}}%
c_{ni}^{2}=1$,  we have for
all $y$, $B_{n}^{2}(-\infty,y)\leq1,$ and therefore, Theorem \ref{FN} implies:%
\[
\mathbb{P}(S_{n}^{(\varepsilon x)}\geq x)\leq\exp(-\alpha^{2}x^{2}%
/2e^{m})+(A_{n}(m;0,\varepsilon x)/(\beta\varepsilon^{m-1}x^{m}))^{\beta
/\varepsilon}\text{ }.
\]
with $\alpha=1-\beta=2/(m+2)$. Then, obviously, it is enough to show that for
$x=x_{n}$ as in Theorem \ref{LargeD1} we can select $\varepsilon>0$ such that
\begin{equation}
\exp(-\frac{\alpha^{2}x^{2}}{2e^{m}})+\left(  \frac{A_{n}(m;0,\varepsilon
x)}{\beta\varepsilon^{m-1}x^{m}}\right)  ^{\beta/\varepsilon}=o(1)\sum
_{i=1}^{k_{n}}\frac{c_{ni}^{t}}{x^{t}}h(\frac{x}{c_{ni}})\text{ as
}n\rightarrow\infty. \label{toshow}%
\end{equation}
Let $x=x_{n}\geq C[\ln(D_{nt}^{-1})]^{1/2}$ where $C>e^{m/2}(m+2)/\sqrt{2}$.
As we mentioned at the beginning of the proof, we clearly have $x_{n}%
\rightarrow\infty.$

We shall estimate each term in the left hand side of (\ref{toshow})
separately. Because, by the definition of $\alpha$ we have $C >e^{m/2}%
\alpha^{-1}\sqrt{2},$ we can select $0<\eta<1$ such that $C ^{2}\alpha
^{2}/2e^{m}=(1-\eta)^{-2}.$

Taking into account the fact that for any $c>0$ and $d>0$ we have $y^{d}%
\exp(-cy)$\newline$=o(\exp(-c(1-\eta)y)$ as $y\rightarrow\infty,$ by the
definition on $x$ and $\eta,$ we obtain:%
\begin{gather*}
x^{(t-2\eta)/(1-\eta)}\exp(-\frac{\alpha^{2}x^{2}}{2e^{m}})=o(1)\exp
(-\frac{\alpha^{2}x^{2}}{2e^{m}}(1-\eta))\\
=o(1){\Large (}\sum_{i=1}^{k_{n}}c_{ni}^{t}{\Large )}^{C^{2}\alpha^{2}%
(1-\eta)/2e^{m}}=o(1){\Large (}\sum_{i=1}^{k_{n}}c_{ni}^{t}{\Large )}%
^{(1-\eta)^{-1}}.
\end{gather*}
 Applying now the H\"{o}lder inequality we
clearly have,
\begin{equation}
\sum_{i=1}^{k_{n}}c_{ni}^{t}=\sum_{i=1}^{k_{n}}c_{ni}^{2\eta}c_{ni}^{t-2\eta
}\leq(\sum_{i=1}^{k_{n}}c_{ni}^{2})^{\eta}(\sum_{i=1}^{k_{n}}c_{ni}%
^{(t-2\eta)/(1-\eta)})^{1-\eta}. \label{rel-cnt}%
\end{equation}
Taking into account that $\sum_{i=1}^{k_{n}}c_{ni}^{2}=1,$ we obtain overall%
\[
\exp(-\frac{\alpha^{2}x^{2}}{2e^{m}})=o(1)x^{-(t-2\eta)/(1-\eta)}\sum
_{i=1}^{k_{n}}c_{ni}^{(t-2\eta)/(1-\eta)}.
\]
Since $t>2,$ $(t-2\eta)/(1-\eta)>$ $t.$ Then, by combining this observation
with the properties of slowly varying functions we have
\[
\exp(-\frac{\alpha^{2}x^{2}}{2e^{m}})=o(1)\sum_{i=1}^{k_{n}}\frac{c_{ni}^{t}%
}{x^{t}}h(\frac{x}{c_{ni}}).
\]
We select $\varepsilon$ by analyzing the second term in the left hand side of
(\ref{toshow}). Notice that by integration by parts formula, for every
$z>y>0,$
\begin{align*}
\mathbb{E}\xi_{0}^{m}I(0  &  \leq\xi_{0}<z)=\\
-z^{m}\mathbb{P}(\xi_{0}  &  \geq z)+m\int_{0}^{z}u^{m-1}\mathbb{P}(\xi
_{0}\geq u)du\leq y^{m}+m\int_{y}^{z}u^{m-1}\mathbb{P}(\xi_{0}\geq u)du.
\end{align*}
Replacing $z=\varepsilon x/c_{ni}$, taking into account condition
(\ref{tail1}), the properties of slowly varying functions, and the facts that
$x/c_{ni}\rightarrow\infty$ and $m>t,$ we have
\[
\mathbb{E}\xi_{0}^{m}I(0\leq c_{ni}\xi_{0}<\varepsilon x)\leq y^{m}+2m\int
_{y}^{\frac{\varepsilon x}{c_{ni}}}u^{m-t-1}h(u)du=O((\frac{x}{c_{ni}}%
)^{m-t}h(\frac{x}{c_{ni}}))
\]
for $y$ sufficiently large. It follows that 
\begin{gather*}
A_{n}(m;0,\varepsilon x)=\sum_{i=1}^{k_{n}}c_{ni}^{m}\mathbb{E}\xi_{0}%
^{m}I(0\leq c_{ni}\xi_{0}<\varepsilon x)\\
\ll\sum_{i=1}^{k_{n}}c_{ni}^{m}(\frac{x}{c_{ni}})^{m-t}h(\frac{x}{c_{ni}})=
x^{m-t}\sum_{i=1}^{k_{n}}c_{ni}^{t}h(\frac{x}{c_{ni}}).
\end{gather*}
Choose $\varepsilon$ with $0<\varepsilon<\beta.$ Then the second term has the
order
\[
\left(  \frac{A_{n}(m;0,\varepsilon x)}{\beta\varepsilon^{m-1}x^{m}}\right)
^{\beta/\varepsilon}\ll\left(  \frac{x^{m-t}}{x^{m}}\sum_{i=1}^{k_{n}}%
c_{ni}^{t}h(\frac{x}{c_{ni}})\right)  ^{\beta/\varepsilon}=o\left(  \sum
_{i=1}^{k_{n}}\frac{c_{ni}^{t}}{x^{t}}h(\frac{x}{c_{ni}})\right)  .
\]

Overall we obtain for any $x\geq C(\ln(\sum_{i=1}^{k_{n}}c_{ni}^{t}%
)^{-1})^{1/2}$ with $C>e^{m/2}(m+2)/\sqrt{2},$%
\[
\mathbb{P}(S_{n}\geq x)=(1+o(1))\sum_{i=1}^{k_{n}}\mathbb{P}(c_{ni}\xi_{0}\geq
x)\text{ as }n\rightarrow\infty,
\]
where $m>t$. Since $C_{t}>e^{t/2}(t+2)/\sqrt{2}$ we can select and fix $m>t$
such that $C_{t}>e^{m/2}(m+2)/\sqrt{2}$. $\lozenge$

\subsection{Proof of Theorem \ref{ModerateD}}

For simplicity we normalize by the variance of $S_{n}$ and assume
(\ref{sum1}). This result easily follows from Theorem 1.1 in Frolov (2005)
when moments strictly larger than $2$ are available. This theorem is given for
convenience in the Appendix (Theorem \ref{frolov}). Because we assume the
existence of moments of order $p>2$, we have
\[
\Lambda_{n}(u,s,\epsilon)\leq u\sum_{j=1}^{k_{n}}c_{nj}^{2}\mathbb{E}\xi
_{0}^{2}I(|c_{nj}\xi_{0}|>\epsilon/s)\leq\epsilon^{2-p}us^{p-2}D_{np}%
\mathbb{E}|\xi_{0}|^{p}.
\]
where $D_{np}=\sum_{j=1}^{k_{n}}|c_{nj}|^{p}$. Then, for $x^{2}\leq
2\ln(1/D_{np}),$ 
\[
\Lambda_{n}(x^{4},x^{5},\epsilon)\leq\epsilon^{2-p}x^{4+5(p-2)}D_{np}%
\mathbb{E}|\xi_{0}|^{p}\leq\epsilon^{2-p}D_{np}(2\ln(1/D_{np}))^{(5p-6)/2}%
\mathbb{E}|\xi_{0}|^{p},
\]
which converges to $0$ since $D_{np}\leq\max_{1\leq j\leq k_{n}}|c_{nj}%
|^{p-2}\rightarrow0$ by (\ref{cni}). Notice also that the $L_{np}$ in Theorem
\ref{frolov} satisfies $L_{np}\leq D_{np}\mathbb{E}|\xi_{0}|^{p}\rightarrow0$.
The latter implies $x^{2}-2\ln(L_{np}^{-1})-(p-1)\ln\ln(L_{np}^{-1}%
)\rightarrow-\infty$ provided $x^{2}\leq2\ln(D_{np}^{-1})$. Then the result is
immediate from Theorem \ref{frolov}. $\lozenge$

\subsection{Proof of Theorem \ref{mix}}

Again for simplicity we normalize by the variance and assume (\ref{sum1}).
Without loss of generality we may assume $2<p<t$. This is so because if $p\ge
t$ with $\mathbb{E(}|\xi_{0}|^{p})<\infty$ then we can find a $p^{\prime}$
such that $2<p^{\prime}<t$ and $\mathbb{E(}|\xi_{0}^{\prime}|^{p})<\infty.$ We
shall consider a sequence $x_{n}$ which converges to $\infty.$ So, let
$x=x_{n}\rightarrow\infty.$

Starting from the relation (\ref{ineqLMD}) and applying
Proposition \ref{frolov-trunc} to the second term in the right hand side we
obtain for any $\varepsilon>0$ and $x^{2}\leq c_{\varepsilon}\ln(D_{np}^{-1})$
with $c_{\varepsilon}<1/\varepsilon$ and for all $n$ sufficiently large
$\mathbb{P}(S_{n}^{(\varepsilon x)}\geq x)=(1-\Phi(x))(1+o(1)).\mathbb{\ }$We
notice now that by (\ref{rel-cnt}) applied with $\eta=(t-p)/(t-2)$ and simple
considerations,%
\begin{equation}
D_{nt}\ll D_{np}\ll(D_{nt})^{(p-2)/(t-2)}. \label{relD}%
\end{equation}
So far, by using this last relation, we showed by (\ref{ineqLMD}) and the
above considerations that (\ref{MD+LD}) holds for $0<x\leq C[\ln(D_{nt}%
^{-1})]^{1/2}$ with $C\ $an arbitrary positive number. On the other hand,
because $1-\Phi(x)\leq (2\pi)^{-1/2}x^{-1}\exp(-x^{2}/2),$ by Theorem \ref{LargeD1} and by
the arguments leading to the proof of relation (\ref{toshow}), there is a
constant $c_{1}>0$ such that for $x>c_{1}[\ln(D_{nt}^{-1})]^{1/2},$ we
simultaneously have
\[
\mathbb{P}\left(  S_{n}\geq x\right)  =(1+o(1))\sum_{i=1}^{k_{n}}%
\mathbb{P}(c_{ni}\xi_{0}\geq x)
\]
and%
\[
1-\Phi(x)=o(\sum_{i=1}^{k_{n}}\mathbb{P}(c_{ni}\xi_{0}\geq x)).
\]
Then (\ref{MD+LD}) holds for all $x>0$ since $C$ is arbitrarily large and can
be selected such that $c_{1}<C$.

Now if the sequence $x_{n}$ is bounded we apply first Theorem \ref{ModerateD}
and obtain the moderate deviation result in (\ref{MD}). Then, because
$x_{n}\geq c>0$ we notice that, by the arguments leading to the proof of
relation (\ref{toshow}), the second part in the right hand side of
(\ref{MD+LD}) is dominant, so the first part is negligible as $n\rightarrow
\infty$. $\lozenge$

\subsection{Proof of Corollary \ref{Remark 1}}

Again without loss of generality we normalize by the variance and assume
(\ref{sum1}). The ideas involved in the proof of this corollary already
appeared in the previous proofs, so we shall mention only the changes. We
start from (\ref{MD+LD}). To prove (\ref{LD}) we have to show that
\[
1-\Phi(x)=o(\sum_{i=1}^{k_{n}}\mathbb{P}(c_{ni}\xi_{0}\geq x))
\]
for $x\geq a(\ln D_{nt}^{-1})^{1/2}$ with $a>2^{1/2}.\ $First we shall use the
relation $1-\Phi(x)\leq (2\pi)^{-1/2}x^{-1}\exp(-x^{2}/2).\ $Then, we adapt the proof we
used to establish the first part of (\ref{toshow}), when we compared
$\exp(-\alpha^{2}x^{2}/2e^{m})$ to $\sum_{i=1}^{k_{n}}\mathbb{P}(c_{ni}\xi
_{0}\geq x).$ The main difference is that now we take $m=0$ and $\alpha=1$.

For the proof of (\ref{MD}), we use the inequality $1-\Phi(x)\geq (2\pi)^{-1/2}
(1+x)^{-1}\exp(-x^{2}/2).$ By (\ref{tail1}) and (\ref{relD}) we have for every
$0<\varepsilon<t-2,$%
\[
\sum_{i=1}^{k_{n}}\mathbb{P}(c_{ni}\xi_{0}\geq x)\ll\sum_{i=1}^{k_{n}}%
\frac{c_{ni}^{t-\varepsilon}}{x^{t-\varepsilon}}\ll\frac{1}{x^{t-\varepsilon}%
}(D_{nt})^{(t-2-\varepsilon)/(t-2)}.
\]
Then, it is easy to see\ that, because $\varepsilon$ can be made arbitrarily
small, for $1<x\leq b(\ln D_{nt}^{-1})^{1/2}$ with $b<2^{1/2}$ we have
\[
\sum_{i=1}^{k_{n}}\mathbb{P}(c_{ni}\xi_{0}\geq x)=o(1-\Phi(x)).
\]
When $0<x\leq1$ we apply Theorem \ref{ModerateD}. $\lozenge$

\subsection{Proof of Corollary \ref{LinearLDMD}}

As in the other proofs, for simplicity we assume $\mathbb{E}\xi_{0}^{2}%
=1.$

\noindent\textbf{Proof of part (ii).} Because the Fuk-Nagaev inequality
(Theorem \ref{FN}) and the inequalities in Lemma \ref{approx} and Proposition
\ref{stockB} are still valid for the case $k_{n}=\infty$ (see Remark
\ref{FNinfty} in the Appendix, Lemma \ref{approxinfty} and Remark
\ref{stochBremark} in Subsection 4.1), all the arguments in the proof of
Theorem \ref{LargeD1}  hold under the conditions
of this corollary. \newline

\noindent\textbf{Proof of part (iii).} The result (iii) in this corollary is
obtained on the same lines as of Theorem \ref{ModerateD}. The modification of
the proof is rather standard but computationally intensive. There are several
ideas behind this proof. The infinite series is decomposed as a sum up to
$k_{n}$ and the rest $R_{n}$. The sequence $k_{n}$ is selected independently
of $x_{n}$ such that the rest of the series $R_{n}$ is negligible for the
moderate deviation result. This is possible because the coefficients $b_{ni},$
defined as $b_{ni}=a_{1-i} +...+a_{n-i}$ with $\sum_{i\in\mathbb{Z}}a_{i}%
^{2}<\infty,$ have some regularity properties. For instance by the
H\"{o}lder's inequality,
\[
b_{ni}^{2}\leq n(a_{1-i}^{2}+...+a_{n-i}^{2})
\]
and so, for any $k>n$
\begin{equation}%
{\displaystyle\sum\nolimits_{|i|\geq k}}
b_{ni}^{2}\leq n^{2}%
{\displaystyle\sum\nolimits_{|i|\geq k-n-1}}
a_{1-i}^{2}. \label{holder}%
\end{equation}

We then note that the existence of moments of order $p>2$ for $\xi_{0}$ and
$\sum_{i\in\mathbb{Z}}a_{i}^{2}<\infty$ imply that $X_{0}$ also has finite
moments of order $p$. Indeed, by Rosenthal inequality (see for instance
Theorem 1.5.13 in de la Pe\~{n}a and Gin\'{e}, 1999), there is a constant
$C_{p}$ such that%
\[
\mathbb{E}|\sum_{j=n}^{m}a_{j}\xi_{j}|^{p}\leq C_{p}[(\sum_{j=n}^{m}a_{j}%
^{2})^{p/2}+\mathbb{E}|\xi_{0}|^{p}\sum_{j=n}^{m}|a_{j}|^{p}]
\]
which implies that $\mathbb{E}|\sum_{j=n}^{m}a_{j}\xi_{j}|^{p}\rightarrow0$ as
$m\geq n\rightarrow\infty,$ and therefore $X_{0}$ exists in $\mathbb{L}_{p}.$

For $k_{n}$ a sequence of integers, denote $R_{n}=\sum_{|i|>k_{n}}b_{ni}%
\xi_{i}$ and note that $R_{n}$ is also well defined in $\mathbb{L}_{p}.$ Again
by Rosenthal inequality we obtain
\begin{equation}
\mathbb{E}|R_{n}|^{p}\leq C_{p}[(\sum_{|i|>k_{n}}b_{ni}^{2})^{p/2}%
+\mathbb{E}|\xi_{0}|^{p}\sum_{|i|>k_{n}}|b_{ni}|^{p}]. \label{ros2}%
\end{equation}
We select now $k_{n}$ large enough such that%
\[
\sum_{|i|>k_{n}}b_{ni}^{2}\leq||\xi_{0}||_{p}^{2}(\sum_{j}|b_{nj}|^{p}%
)^{2/p}.
\]
This is possible by relation (\ref{holder}) and the fact that $\sum_{i}%
a_{i}^{2}<\infty$. With this selection we obtain%
\begin{equation}
\mathbb{E}|R_{n}|^{p}\leq2C_{p}\mathbb{E}|\xi_{0}|^{p}\sum_{i}|b_{ni}|^{p}.
\label{p-rest}%
\end{equation}
Write now
\[
S_{n}=\sum_{|i|\leq k_{n}}b_{ni}\xi_{i}+R_{n}.
\]
We view $S_{n}$ as the sum of $k_{n}+1$ independent random variables and then
apply Theorem \ref{frolov} as in the proof of Theorem \ref{ModerateD}. By
taking into account (\ref{p-rest}), the term $L_{np}$ from Theorem
\ref{frolov} is 
\begin{align*}
L_{np}  &  =\frac{1}{\sigma_{n}^{p}}[\sum_{|i|\leq k_{n}}|b_{ni}%
|^{p}\mathbb{E(}\xi_{0}^{p}I(\xi_{0}>0)+\mathbb{E(}R_{n}^{p}I(R_{n}>0)]\\
&  \leq\frac{2C_{p}+1}{\sigma_{n}^{p}}\sum_{i}|b_{ni}|^{p}\mathbb{E}|\xi
_{0}|^{p}=(2C_{p}+1)U_{np}\mathbb{E}|\xi_{0}|^{p}=L_{np}^{\prime}.
\end{align*}
Because we assume the existence of moments of order $p$, by (\ref{p-rest}) we
have
\begin{gather*}
\Lambda_{n}(u,s,\epsilon)\leq\frac{u}{\sigma_{n}^{2}}[\sum_{|j|\leq k_{n}%
}b_{nj}^{2}\mathbb{E}\xi_{0}^{2}I(|b_{nj}\xi_{0}|>\epsilon\sigma
_{n}/s)+\mathbb{E}R_{n}^{2}I(|R_{n}|>\epsilon\sigma_{n}/s)]\\
\leq\frac{us^{p-2}}{\sigma_{n}^{p}\epsilon^{p-2}}[\sum_{|j|\leq k_{n}}%
|b_{nj}|^{p}\mathbb{E}|\xi_{0}|^{p}+\mathbb{E}|R_{n}|^{p}]\leq\frac{us^{p-2}%
}{\epsilon^{p-2}}L_{np}^{\prime}.
\end{gather*}
Therefore, for $x^{2}\leq2\ln(1/L_{np}^{\prime})\leq2\ln(1/L_{np}),$
\[
\Lambda_{n}(x^{4},x^{5},\epsilon)\leq\epsilon^{2-p}x^{4+5(p-2)}L_{np}^{\prime
}\mathbb{\ }\leq\epsilon^{2-p}(2\ln(1/L_{np}^{\prime}))^{(5p-6)/2}%
L_{np}^{\prime}\text{.}%
\]
Finally note that by (\ref{bni}) we obtain
\[
U_{np}\leq\frac{\sup_{j}|b_{nj}|^{p-2}}{(\sum b_{nj}^{2})^{(p-2)/2}%
}\rightarrow0\text{,}%
\]
and consequently $L_{np}^{\prime}\rightarrow0.$ Therefore, $\Lambda_{n}%
(x^{4},x^{5},\epsilon)\rightarrow0$. Note also that the quantity $L_{np}$ in
Theorem \ref{frolov} satisfies $L_{np}\leq L_{np}^{\prime}\rightarrow0$.
Therefore if $x^{2}-2\ln(L_{np}^{\prime})^{-1}-(p-1)\ln\ln(L_{np}^{\prime
})^{-1}\rightarrow-\infty$ we have that $x^{2}-2\ln(L_{np}^{-1})-(p-1)\ln
\ln(L_{np}^{-1})\rightarrow-\infty$ and the result holds for such a positive
$x$.

It remains to show that $x^{2}\leq2\ln(U_{np}^{-1})$ implies $x^{2}%
-2\ln(L_{np}^{\prime})^{-1}-(p-1)\ln\ln(L_{np}^{\prime})^{-1}\rightarrow
-\infty,$ which holds provided that%
\[
2\ln((U_{np}^{-1})L_{np}^{\prime}[\ln(L_{np}^{\prime})^{-1}]^{(1-p)/2}%
)\rightarrow-\infty.
\]
This last divergence is equivalent to
\[
(U_{np}^{-1})L_{np}^{\prime}[\ln(L_{np}^{\prime})^{-1}]^{(1-p)/2}%
\rightarrow0.
\]
Clearly, because $L_{np}^{\prime}=(2C_{p}+1)U_{np}\mathbb{E}|\xi_{0}|^{p}$
 and the fact that we have shown that
$L_{np}^{\prime}\rightarrow0$ the result follows. \newline

\noindent\textbf{Proof of part (i).} The proof is similar to the proof of
Theorem \ref{mix} and Corollary \ref{Remark 1}. We have only to show that Proposition \ref{frolov-trunc} is
still valid in this context if we let $k_{n}=\infty$. The proof is similar to
the proof of (iii) but more involved, since the sequence of truncated
variables is not centered. Denote
\[
X_{ni}^{\prime}=b_{ni}\xi_{i}I(b_{ni}\xi_{i}\leq\varepsilon x\sigma
_{n})=b_{ni}\xi_{i}^{\prime}.
\]
For $k_{n}$ a sequence of integers, denote $R_{n}^{\prime}=\sum_{|i|>k_{n}%
}b_{ni}\xi_{i}^{\prime}$ and note that $R_{n}^{\prime}$ is also well defined
in $\mathbb{L}_{p}.$ By Rosenthal inequality, after centering we obtain
\[
\mathbb{E}|R_{n}^{\prime}|^{p}\leq C_{p}^{\prime}[(\sum_{|i|>k_{n}}b_{ni}%
^{2})^{p/2}+\mathbb{E}|\xi_{0}|^{p}\sum_{|i|>k_{n}}|b_{ni}|^{p}+|\mathbb{E}%
(R_{n}^{\prime})|^{p}].
\]
Because $x\geq c>0$ and the fact that $\mathbb{E}(X_{ni}^{\prime}%
)=-\mathbb{E}(b_{ni}\xi_{i}I(b_{ni}\xi_{i}>\varepsilon x\sigma_{n}))$ we
obtain
\[
|\mathbb{E}(R_{n}^{\prime})|\leq\frac{1}{\varepsilon x\sigma_{n}}%
\sum_{|i|>k_{n}}b_{ni}^{2}\leq\frac{1}{\varepsilon c\sigma_{n}}\sum
_{|i|>k_{n}}b_{ni}^{2}.
\]
We select now $k_{n}$, depending on $c,$ $\varepsilon$ and the distribution of
$\xi_{0}$ and the coefficients $(a_{k}),$ large enough such that%

\[
(\sum_{|i|>k_{n}}b_{ni}^{2})^{p/2}+(\frac{1}{\varepsilon c\sigma_{n}}%
\sum_{|i|>k_{n}}b_{ni}^{2})^{p}\leq\mathbb{E}|\xi_{0}|^{p}\sum_{i}|b_{ni}%
|^{p},
\]
and so
\[
\mathbb{E}|R_{n}^{\prime}|^{p}\leq2C_{p}^{\prime}\mathbb{E}|\xi_{0}|^{p}%
\sum_{i}|b_{ni}|^{p}.
\]

Write now $S_{n}^{\prime}=\sum_{|i|\leq k_{n}}b_{ni}\xi_{i}^{\prime}%
+R_{n}^{\prime}$ and view $S_{n}^{\prime}$ as the sum of $k_{n}+1$ independent
random variables and then apply Proposition \ref{frolov-trunc}. Similar
computations as in the proof of the point (iii) show that $L_{np}$ in
Proposition \ref{frolov-trunc} is bounded by 
\[
L_{np}\leq\frac{2C_{p}^{\prime}+1}{\sigma_{n}^{p}}\mathbb{E}|\xi_{0}|^{p}%
\sum_{i}|b_{ni}|^{p}=(2C_{p}^{\prime}+1)U_{np}\mathbb{E}|\xi_{0}|^{p}.
\]
Then, by Proposition \ref{frolov-trunc} if $x^{2}\leq c\ln((2C_{p}^{\prime
}+1)U_{np}\mathbb{E}|\xi_{0}|^{p})^{-1}$ for $c<1/\varepsilon$, we have
$x^{2}\leq c\ln(L_{np}^{-1})$ for $c<1/\varepsilon$ and
\[
\mathbb{P}\left(
{\displaystyle\sum\nolimits_{i}}
X_{nj}^{\prime}\geq x\sigma_{n}\right)  =(1-\Phi(x))(1+o(1)).
\]
It remains to notice that because $U_{np}\rightarrow0,$ we also have the
result for $x^{2}\leq c\ln(U_{np})^{-1}$ for any $c<1/\varepsilon$, for all
$n$ sufficiently large. $\lozenge$

\subsection{Proof of Corollary \ref{RegularLDMD}}

This Corollary follows from Corollary \ref{LinearLDMD} via Lemma \ref{coeff}
in the Appendix. It remains to give an explicit form of the intervals moderate
deviation and large deviation boundaries. Without loss of generality, we
assume that $\mathbb{E}\xi_{0}^{2}=1$. For proving the large deviation part of
this corollary we have to analyze the condition on $x$ from part (i) of
Corollary \ref{LinearLDMD}, namely $x>a(\ln U_{nt}^{-1})^{1/2}$ with
$a=\sqrt{2}$. By Lemma \ref{coeff}
\[
B_{n2}=\sum_{i}b_{ni}^{2}\sim c_{r}n^{3-2r}l^{2}(n)
\]
and
\[
C_{1}l^{t}(n)n^{(1-r)t+1}\leq\sum_{j=1}^{\infty}b_{nj}^{t}\leq C_{2}%
l^{t}(n)n^{(1-r)t+1}\text{.}%
\]
Then, for certain constants $K_{1}$ and $K_{2}$ and because $U_{nt}%
^{-1}=B_{n2}^{t/2}/B_{nt},$ we have for $n$ sufficiently large
\[
K_{1}+\ln n^{(t-2)/2}\leq\ln U_{nt}^{-1}\leq K_{2}+\ln n^{(t-2)/2}.\text{ }%
\]
So, the asymptotic result (\ref{LD}) holds for $x\geq c_{1}(\ln n)^{1/2}$
where $c_{1}>(t-2)^{1/2}.$ Furthermore, (\ref{MD}) holds for $0<x\leq
c_{2}(\ln n)^{1/2}$ where $c_{2}<(t-2)^{1/2}$. $\lozenge$

\subsection{Proof of Theorem \ref{thmfunc}}

\label{sec:prf24} Without restricting the generality we assume $\kappa>0,$
since similar computations can be done when $\kappa<0.$ Let $A_{n}=\sum
_{i=n}^{\infty}a_{i}^{2}$. Using the argument of Theorem 5 in Wu (2006), under
Condition B, we have
\[
\Vert\mathcal{P}_{0}(K(X_{n})-\kappa X_{n})\Vert_{q}=O(\theta_{n}),\text{
where }\theta_{n}=|a_{n}|^{p/q}+|a_{n}|A_{n}^{1/2}.
\]
Let $\theta_{i}=0$ if $i\leq0$ and $\Theta_{n}=\sum_{i=1}^{n}\theta_{i}$. Then
by Theorem 1 in Wu (2007), there exists a constant $B_{q}\geq1$ such that
\begin{equation}
{\frac{{\Vert S_{n,1}\Vert_{q}^{2}}}{{B_{q}^{2}}}}\leq\sum_{i\in\mathbb{Z}%
}(\Theta_{n+i}-\Theta_{i})^{2}\leq2n\Theta_{2n}^{2}+\sum_{i=n+1}^{\infty
}(\Theta_{n+i}-\Theta_{i})^{2}. \label{momentq}%
\end{equation}
By Karamata's theorem, $A_{n}\sim(2r-1)^{-1}n^{1-2r}l(n)^{2}$, and if $i>n$,
$\Theta_{n+i}-\Theta_{i}=O(n\theta_{i})$ and $\sum_{i=n+1}^{\infty}\theta
_{i}^{2}=O(n\theta_{n}^{2})$. Let $\ell(\cdot)$ be a slowly varying function
and $\beta\in\mathbb{R}$. Again by Karamata's theorem, there exists another
slowly varying function $\ell_{0}(\cdot)$ such that $\sum_{i=1}^{n}i^{-\beta
}\ell(i)=O(1+n^{1-\beta})\ell_{0}(n)$. Hence by (\ref{momentq}), there exists
a slowly varying function $\ell_{1}(\cdot)$ such that
\begin{equation}
\Vert S_{n,1}\Vert_{q}=O(\sqrt{n})(1+n^{1-rp/q}+n^{1-r+(1-2r)/2})\ell_{1}(n).
\label{16}%
\end{equation}
For $n\geq3$ let $g_{n}=(\ln n)^{-1}$. Then
\begin{equation}
\mathbb{P}(S_{n}\geq(x+g_{n})\sigma_{n})-\mathbb{P}(H_{n}\geq\kappa
x\sigma_{n})\leq\mathbb{P}(|S_{n,1}|\geq\kappa g_{n}\sigma_{n}). \label{diff}%
\end{equation}
Since $x^{2}\leq c\ln n$ and $g_{n}=(\ln n)^{-1}$, we have that $1-\Phi(x\pm
g_{n})\sim1-\Phi(x)$. Hence by Corollary \ref{Separate}, (\ref{MDH}) follows
from (\ref{diff}) in view of
\begin{align}
\mathbb{P}(|S_{n,1}|  &  \geq\kappa g_{n}\sigma_{n})\leq{\frac{{\Vert
S_{n,1}\Vert_{q}^{q}}}{|\kappa|^{q}{g_{n}^{q}\sigma_{n}^{q}}}}={\frac
{{O(\sqrt{n}^{q})(1+n^{q-rp}+n^{(3/2-2r)q})\ell_{1}^{q}(n)}}{{g_{n}%
^{q}(n^{3/2-r}l(n))^{q}}}}\label{18}\\
&  =n^{-p\rho(r)}{\frac{{\ell_{1}^{q}(n)}}{{g_{n}^{q}l}^{q}{(n)}}}%
={\frac{{o(n^{-c/2})}}{\ln n}}=o(xe^{-x/2})=o[1-\Phi(x)],\nonumber
\end{align}
since $c/2<p\rho(r)$. Here we note that $\ell_{1}(n)/(g_{n}l(n))$ is also
slowly varying in $n$ and $x\leq c\ln n$. By (\ref{16}) and (\ref{18}), it is
easily seen that the normalizing constant $\kappa\sigma_{n}$ can be replaced
by $\sqrt{var(H_{n})}$. The proof of the upper bound is similar and it is left
to the reader. $\lozenge$

\section{Appendix}

The following Theorem is a slight reformulation of Fuk--Nagaev inequality (see
Corollary 1.7, S. Nagaev, 1979):

\begin{theorem}
\label{FN}Let $X_{1},\cdots,X_{k_{n}}$ be independent random variables. Assume
$m\geq2$. Suppose $\mathbb{E}X_{i}=0,$ $i=1,\cdots,k_{n},$ $\beta=m/(m+2)$,
and $\alpha=1-\beta=2/(m+2)$. For $y>0$, define $X_{i}^{(y)}=X_{i}I(X_{i}\leq
y)$, $A_{n}(m;0,y):=\sum_{i=1}^{k_{n}}\mathbb{E}[X_{i}^{m}I(0<X_{i}<y)]$ and
$B_{n}^{2}(-\infty,y):=\sum_{i=1}^{k_{n}}\mathbb{E}[X_{i}^{2}I(X_{i}<y)].$
Then for any $x>0$ and $y>0$
\begin{equation}
\mathbb{P}(\sum_{i=1}^{k_{n}}X_{i}^{(y)}\geq x)\leq\exp(-\frac{\alpha^{2}%
x^{2}}{2e^{m}B_{n}^{2}(-\infty,y)})+{\LARGE (}\frac{A_{n}(m;0,y)}{\beta
xy^{m-1}}{\LARGE )}^{\beta x/y}. \label{Nagaev2}%
\end{equation}

\end{theorem}

\begin{remark}
\label{FNinfty} Let $X_{1},X_{2},\cdots,$ be independent random variables.
Assume that the sum $S=\sum_{i=1}^{\infty}X_{i}$ exists almost surely. By the
same argument as in Lemma \ref{approxinfty}, $\sum_{i=1}^{\infty}X_{i}^{(y)}$
converges almost surely for all $y>0$. By passing to the limit in
(\ref{Nagaev2}) we note that this version of Fuk-Nagaev inequality is still
valid for $\mathbb{P}(\sum_{i=1}^{\infty}X_{i}^{(y)}\geq x)$.
\end{remark}

We shall also use the following result which is an immediate consequence of
Theorem 1.1 in Frolov (2005).

\begin{theorem}
\label{frolov} Let $(X_{nj})_{1\leq j\leq k_{n}}$ be an array of row-wise
independent centered random variables. Let $p>2$ and denote $S_{n}=\sum
_{j=1}^{k_{n}}X_{nj}$, $\sigma_n^2=\sum_{j=1}^{k_n}\mathbb{E}X_{nj}^2$, $M_{np}=\sum_{j=1}^{k_{n}}\mathbb{E}X_{nj}^{p}%
I(X_{nj}\geq0)<\infty$, $L_{np}=\sigma_{n}^{-p}M_{np}$ and denote
\[
\Lambda_{n}(u,s,\epsilon)=\frac{u}{\sigma_{n}^{2}}\sum_{j=1}^{k_{n}}%
\mathbb{E}X_{nj}^{2}I(X_{nj}\leq-\epsilon\sigma_{n}/s).
\]
Furthermore, assume $L_{np}\rightarrow0$ and $\Lambda_{n}(x^{4},x^{5}%
,\epsilon)\rightarrow0$ for any $\epsilon>0$. Then if $x\geq0$ and $x^{2}%
-2\ln(L_{np}^{-1})-(p-1)\ln\ln(L_{np}^{-1})\rightarrow-\infty,$ we have
\[
\mathbb{P}\left(  S_{n}\geq x\sigma_{n}\right)  =(1-\Phi(x))(1+o(1)).
\]

\end{theorem}

For truncated random variables by following the proof of Theorem 1.1 in Frolov
(2005) we can present his relation (3.17) as a proposition.

\begin{proposition}
\label{frolov-trunc} Assume the conditions in Theorem \ref{frolov} are
satisfied. Fix $\varepsilon>0.$ Define
\[
X_{nj}^{(\varepsilon x\sigma_{n})}=X_{nj}I(X_{nj}\leq\varepsilon x\sigma
_{n})\text{ and }S_{n}^{(\varepsilon x\sigma_{n})}=\sum_{j=1}^{k_{n}}%
X_{nj}^{(\varepsilon x\sigma_{n})}.
\]
Then if $x^{2}\leq c\ln(L_{np}^{-1})$ with $c<1/\varepsilon$, for all $n$
sufficiently large we have
\[
\mathbb{P}\left(  S_{n}^{(\varepsilon x\sigma_{n})}\geq x\sigma_{n}\right)
=(1-\Phi(x))(1+o(1)).
\]

\end{proposition}


The following facts about the series are going to be used to analyze a class
of linear processes:

\begin{lemma}
\label{coeff}Assume $a_{i}=l(i)i^{-r}$ with $1/2<r<1$. Let $b_{j}%
:=b_{nj}:=\sum_{i=1}^{j}a_{i}$ if $1\leq j\leq n$ and $b_{nj}:=\sum
_{i=j-n+1}^{j}a_{i}$ if $j>n$. Then, for two positive constants $C_{1}$ and
$C_{2},$ we have
\[
C_{1}(l^{t}(n)n^{(1-r)t+1})\leq\sum_{j=1}^{\infty}b_{nj}^{t}\leq C_{2}%
(l^{t}(n)n^{(1-r)t+1}),
\]
for any $t\geq2$. In the case $t=2$, $\sum_{j=1}^{\infty}b_{nj}^{2}%
=c_{r}n^{3-2r}l^{2}(n)$ with
\[
c_{r}=\{\int_{0}^{\infty}[x^{1-r}-\max(x-1,0)^{1-r}]^{2}dx\}/(1-r)^{2}.
\]

\end{lemma}

\textbf{Proof.} It is easy to see that $b_{nj}\ll j^{1-r}l(j)$ for $j\leq2n$
and $b_{nj}\ll n(j-n)^{-r}l(j)$ for $j>2n$ from the Karamata theorem (see part
1 of Lemma 5.4 in Peligrad and Sang (2012)). Therefore,
\begin{align*}
&  \sum_{j=1}^{\infty}b_{nj}^{t}=\sum_{j=1}^{2n}b_{nj}^{t}+\sum_{j=2n+1}%
^{\infty}b_{nj}^{t}\\
&  \ll\sum_{j=1}^{2n}j^{(1-r)t}l^{t}(j)+\sum_{j=2n+1}^{\infty}n^{t}%
(j-n)^{-rt}l^{t}(j)=O(l^{t}(n)n^{(1-r)t+1}).
\end{align*}
The proof in the other direction is similar. The result of case $t=2$ is well
known. See for instance Theorem 2 in Wu and Min (2005). $\lozenge$

\bigskip

\textbf{Acknowledgement.} The authors would like to thank the referees for
carefully reading the manuscript and for numerous suggestions that improved
the presentation.

\end{document}